\documentclass[12pt, notitlepage]{article}
\usepackage[utf8]{inputenc}
\usepackage[english]{babel}
\usepackage{amsfonts}
\usepackage{amssymb}
\usepackage{amsmath}
\usepackage{amsthm}
\usepackage{tocloft}
\usepackage{CJK}
\usepackage{tikz-cd}
\usepackage{mathrsfs}
\usepackage[toc,page]{appendix}

\usepackage[english]{babel}

\usepackage{bm} 
\usepackage{chngcntr}
\makeatletter
\def\blfootnote{\xdef\@thefnmark{}\@footnotetext}
\makeatother

\counterwithin{equation}{section}

\makeatletter\@addtoreset{chapter}{part}\makeatother

\newcommand{\xdownarrow}[1]{%
  {\left\downarrow\vbox to #1{}\right.\kern-\nulldelimiterspace}
}

\begin{document}

\title{Real intersection theory I}

\author{B. Wang\\
\begin{CJK}{UTF8}{gbsn}
(汪      镔)
\end{CJK}}

\maketitle

\begin{abstract}
We introduce the notion of  Lebesgue currents. 
They are a special type of currents involving Lebesgue measure. 
We apply it to 
define the intersection of singular cycles, which serves as  the foundation for the 
real intersection theory. 

\end{abstract}

\newcommand{\hookuparrow}{\mathrel{\rotatebox[origin=c]{90}{$\hookrightarrow$}}}
\newcommand{\hookdownarrow}{\mathrel{\rotatebox[origin=c]{-90}{$\hookrightarrow$}}}
\newcommand\ddaaux{\rotatebox[origin=c]{-90}{\scalebox{0.70}{$\dashrightarrow$}}} 
\newcommand\dashdownarrow{\mathrel{\text{\ddaaux}}}

\newtheorem{thm}{Theorem}[section]

\newtheorem{ass}[thm]{\bf {Claim} }
\newtheorem{prop}[thm]{\bf {Property } }
\newtheorem{prodef}[thm]{\bf {Proposition and Definition } }
\newtheorem{construction}[thm]{\bf {Main construction } }
\newtheorem{assumption}[thm]{\bf {Assumption} }
\newtheorem{proposition}[thm]{\bf {Proposition} }
\newtheorem{theorem}[thm]{\bf {Theorem} }
\newtheorem{apd}[thm]{\bf {Algebraic Poincar\'e duality} }
\newtheorem{cond}[thm]{\bf {Condition} }
\newtheorem{ex}[thm]{\bf Example}
\newtheorem{corollary}[thm]{\bf Corollary}
\newtheorem{definition}[thm]{\bf Definition}
\newtheorem{lemma}[thm]{\bf Lemma}
\newtheorem{con}[thm]{C}
\newtheorem{conj}[thm]{\bf Conjecture}
\counterwithin{equation}{section}

\bigskip

\tableofcontents


\section{Introduction}  The idea is originated  from  an attempt  to define a novel   intersection of  singular cycles,  that extends the intersection of algebraic cycles and has a descend to the product in cohomology and Chow groups. 
However as it progresses, it turns away from the original field and emerges  as  the real analysis in geometric measure theory
 which is the content of this paper.      \bigskip

\subsection { History} \par
   Let $\mathcal X$ be a  $C^\infty$ manifold of dimension $m$ (all manifolds in this paper are connected and oriented, but not necessarily compact).
In [3], G. de Rham defined the intersection
\begin{equation}
T\wedge \omega
\end{equation}
between a  current $T$ and  a $C^\infty$ form $\omega$, expressed as a functional on $\mathscr D(\mathcal X)$ --
the space of $C^\infty$ differential forms
with a compact support,  
\begin{equation}
\int_T \omega\wedge \bullet 
\end{equation}
where the integral notation $\int_{T}(\bullet)$ denotes the functional.
The intersection satisfies
\begin{equation}
supp(T\wedge \omega)\subset supp(T)\cap supp(\omega)
\end{equation}where $supp(\cdot)$ stands for the support.   The  asymmetric expression (1.1) calls for a completion that historically 
was dragged into the topology.  For instance, 
 G. de Rham  extended (1.1) to the intersection between two currents, called Kronecker index, 
\begin{equation}
T\wedge S[1],
\end{equation}
where $S$ is another current of dimension $m-dim(T)$ and one of currents has a compact support. 
   To do that, he  first constructed  the De Rham's regularization
$R_\epsilon T$, which is a family of $C^\infty$ forms for real $\epsilon>0$, that converges to
$T$ as $\epsilon\to 0$.   Then he studied the convergence of the real number,
\begin{equation}
\int_{\mathcal X} R_\epsilon(T)\wedge R_{\epsilon'}(S), \quad as \ (\epsilon, \epsilon')\to (0, 0).
\end{equation}
Today we know that the convergence 
 leads  to  the cap product on homology.    
\bigskip

\subsection{ New result} \par
In this paper we  study  the convergence of a similar, but different real number, 
\begin{equation}
\int_{T_1} R_\epsilon T_2\wedge \phi, \quad\quad as \ \epsilon\to 0
\end{equation}
for a fixed test form $\phi\in \mathscr D(\mathcal X)$, where $T_1, T_2$ are currents satisfying
$$dim(T_1)+dim(T_2)\geq m.$$ The convergence would've been a generalization of (1.5), and the limit  could've been written as 
\begin{equation}
T_1\wedge T_2[\phi]
\end{equation}
as the symmetric completion of (1.1).   
Nevertheless the limit of (1.6) does not always exist. 
  To achieve the convergence,  we must avoid the the divergence occurring in the support of currents.
There are multiple directions in doing so, 
 one of which, mentioned above, is De Rham's homotopy approach.  Historically this direction  led to 
his celebrated theorem in differential topology.  But we focus on another direction:  the convergence of (1.6) in analysis yields 
a novel  intersection where the geometry is not required. 
But the  static convergence 
is limited to a special type of currents 
called Lebesgue currents \footnote{ The convergence is intrinsic as an invariant of differential structures, but the limit of the convergence is not.}.
  The tool  is the geometric measure theory. We are trying to see what measure on the support of currents results in the convergence of (1.6). 
In technique, our convergence uses measure theory to
  focus on the growth of measures, contrary to  other known cases where the focus is the growth
 of integrands (i.e. the functions. For instance see [2]). 
Applying the resulting convergence we obtain 
 an extrinsic bilinear homomorphism,  
\begin{equation}\begin{array}{ccc}
\mathcal C(\mathcal X)\times \mathcal C(\mathcal X) &\rightarrow &  \mathcal C(\mathcal X) \\
(T_1, T_2) &\rightarrow & [T_1\wedge T_2],
\end{array}\end{equation}
on the subspace  $\mathcal C(\mathcal X)$ consisting of  Lebesgue currents. So (1.8) does not only  extend the formula (1.1),
but also (1.3), i.e. 
\begin{equation}
supp( [T_1\wedge T_2])\subset supp(T_1)\cap supp(T_2).
\end{equation}
 The property (1.9) is  the original meaning of  
intersection of singular cycles.  
Thus  our  intersection of singular cycles   is a Lebesgue current, which indicates that the cohomoogy ring defined in algebra is related to the geometric measure in the local charts.  Coupled with the homotopy formula associated with  De Rham's regularization, the convergence will lead to a structure beyond De Rham. 
\medskip

In this paper let us  concentrate on the technical part of the convergence of (1.6).  

\bigskip
 
{\bf Acknowledgment} We thank P. Deligne for pointing out the missteps   in the original version of the paper.

\bigskip

\section{Lebesgue currents}
\bigskip

 \bigskip

\subsection{Definition}

\begin{definition} ( Radon-Nikodym\footnote {Radon-Nikodym derivative is an important locally $L^1$ function in the theory of probability (see [1]),   whose  average values around the  non a.e. points lie in the heart of
convergence of (1.6).  } )  In the following,  vectors or points in Euclidean space
$\mathbb R^{\bullet}$ 
will be denoted by ${\bf bold}$ letters. 
Let $\mathbb R^r$ be the Euclidean space with the standard linear structure that has a  
 basis $\mathbf e_{1}, \cdots, \mathbf e_{r}$ and  coordinates $x=\{x_1, \cdots, x_r\}$. 
 Let
$d\mu_x$ be the Euclidean volume form
$$dx_1\wedge\cdots \wedge dx_r, $$

\begin{equation} \bm \lambda=\lambda_{1}\mathbf e_{1}+\cdots+\lambda_{r}\mathbf e_{r}\in \mathbb R^r\end{equation}
be a varied vector in the open region such that 
$\lambda_i>0, all\ i$.  We'll describe a ``path" convergence of a function of $\bm\lambda$ as $\bm\lambda\to \mathbf 0$ in the following way. 
First we divide the basis $\mathbf e_{1}, \cdots, \mathbf e_{r}$ into $l$ groups 
 with an order:  $j_1$ group , $j_2$ group, $\cdots$,  $j_l$ group.  It will be referred to as 
a \begin{equation}
group \quad order
\end{equation}
 for $\mathbb R^r$. 
Then 
we consider such $\bm\lambda$ that  coefficients $\lambda_i$ of $\mathbf e_i$ in the same group are equal to the same number 
$l_{j_k}\in\mathbb R^+ $.  Thus the group order is just an order
 $$\lambda_{j_1}, \lambda_{j_2}, \cdots, \lambda_{j_l}$$
for $\lambda_{j_i}$ near $0$. 
We'll use the symbol  $\lim_{\bm \lambda\Rsh  \mathbf 0}$ to denote the  
 particularly ordered limit for $\bm\lambda \to \mathbf 0$ (i.e. all $\lambda_{j_k}\to 0$) in the order 
$$\displaystyle {\lim_{\lambda_{j_l}\to 0}}\cdots \displaystyle {\lim_{\lambda_{j_1}\to 0}}.$$
  We called the limit, a zigzag limit.
 Let $$\mathbf u=u_{1}\mathbf e_{1}+\cdots+u_{r}\mathbf e_{r}\in \mathbb R^r$$ be a point.
Let $D_{\bm \lambda}$ be an invertible  affine map in the form,
\begin{equation}\begin{array}{ccc}
\mathbb R^r &\rightarrow &\mathbb R^r\\
\mathbf x & \rightarrow  & B\circ  \mathbb D_{\bm\lambda}(\mathbf x)+\mathbf u
 \end{array}\end{equation}
referred to as  the testing map, where $B$ is an invertible linear map and 
$ \mathbb D_{\bm\lambda}$ is the diagonal linear map
$$\begin{array}{ccc}
\mathbb R^r &\rightarrow &\mathbb R^r\\
\mathbf e_i & \rightarrow  & \lambda_i \mathbf e_i, all\ i.
 \end{array}$$

 We say $L\in \mathscr L^1_{loc}(\mathbb R^r)$ is of 
Radon-Nikodym  if 
for any test function $\phi \in \mathscr D(\mathbb R^r)$, any testing map $D_{\bm\lambda}$ and any group order of the basis, 
  the zigzag limit
\begin{equation}\begin{array}{c}
\displaystyle{\lim_{\bm \lambda\Rsh  \mathbf 0}} \int_{\mathbf x\in \mathbb R^r} 
 L ( D_{\bm \lambda }(\mathbf x) )\phi(\mathbf x) d\mu_x\end{array}\end{equation}
exists.   We denote the number (2.4)  by 
\begin{equation}
RN_{\phi, L},
\end{equation}
and call it the Radon-Nikodym number.

\end{definition}

\bigskip

\bigskip

\begin{proposition}
 It does not depend on coordinates for $L\in \mathscr L_{loc}^1(\mathbb R^r)$ to be of  Radon-Nikodym.

\end{proposition}

\bigskip

{\bf Remark} However the Radon-Nikodym number depends the  coordinates, $B$ and the zigzag limiting path $\Rsh$ (i.e.  the group
order). 
\bigskip

\begin{proof} 
\bigskip 

The  proof for $\mathbf u\neq 0$ is identical  with the homogeneous case where $\mathbf u=0$. 
So let's prove the homogeneous case.  Let $L\in \mathscr L_{loc}^1(\mathbb R^r)$ be of 
Radon-Nikodym in $x$-coordinates. 
Let  $y=\{y_1, \cdots, y_r\}$ be another coordinates of $\mathbb R^r$, and 
\begin{equation}\begin{array}{ccc}
\nu: \mathbb R^r &\rightarrow & \mathbb R^r\\
(x_1, \cdots, x_r) &\rightarrow & (y_1, \cdots, y_r) 
\end{array}\end{equation}
be the diffeomorphism between the $x$-$y$ coordinates. 
So we  assumed the homogeneous case, $$\nu(\mathbf 0)=\mathbf 0.$$ 

 Denote the 
volume forms of $\mathbb R^r$ in $y$, $x$ coordinates by $d\mu_y, d\mu_x$ respectively. So 
$$d\mu_y=g(\mathbf x) d\mu_x,$$ 
where $g(\mathbf x)$ is $C^\infty$. 
Then the composition $L\circ \nu^{-1}$ denoted by $L_y$ is  also locally $L^1$.
It is sufficient  to show the convergence of  the number
\begin{equation}
A_{\bm\lambda}=\int_{\mathbf y\in \mathbb R^r}  L_y (D_{\bm\lambda}(\mathbf y))\phi(\mathbf y) d\mu_y
\end{equation}
as $\bm\lambda\Rsh \mathbf 0$, where $D_{\bm \lambda}$ is the testing  map with the linear transformation $B$ and
$\phi(\mathbf y)\in \mathscr D(\mathbb R^r)$.  
First we use standard calculation to convert the expression to $x$-coordinates,
\begin{align}
A_{\bm\lambda}
& ={1\over  det( B )\prod_{i=1}^r \lambda_i}
 \int_{\mathbf y\in \mathbb R^r} L_y(\mathbf y) \phi (D_{\bm\lambda}^{-1} (\mathbf y))  d\mu_y\\
&
={1\over det( B )\prod_{i=1}^r \lambda_i}  \int_{\mathbf x\in \mathbb R^r}  L_y(\nu(\mathbf x)) 
\biggl( \nu^\ast (\phi (D_{\bm\lambda}^{-1} (\mathbf y))  d\mu_y)\biggr)\\ &
={1\over det( B )\prod_{i=1}^r \lambda_i} \int_{\mathbf x\in \mathbb R^r}  L(\mathbf x) 
\biggl( \nu^\ast (\phi (D_{\bm\lambda}^{-1} (\mathbf y))  d\mu_y)\biggr)
\end{align}

We make a change of variable
$$\mathbf x\Rightarrow\alpha\circ D_{\bm\lambda}(\mathbf x)$$
(replacement of $\mathbf x$ with $\alpha\circ D_{\bm\lambda}(\mathbf x)$)
where $\alpha= \nu^{-1}_\ast|_{\mathbf 0}$, a constant matrix under the $y$-basis.  So $det(\alpha) g(\mathbf 0)=1$.
The  integral in the last row  (2.10) is
\begin{equation}
 det(\alpha)\int_{\mathbf x\in \mathbb R^r} L \bigl(\alpha\circ D_{\bm\lambda}  (\mathbf x)\bigr) 
\phi \bigl(D_{\bm\lambda}^{-1}\circ \nu\circ\alpha\circ D_{\bm\lambda}  (\mathbf x) \bigr)
g\bigl(\alpha\circ D_{\bm\lambda}  (\mathbf x)\bigr)
  d\mu_{x}
\end{equation}

Because $\phi$ has a compact support, the variable $\mathbf x$ in the integral (2.11) is bounded.
Hence as $|\bm\lambda|\to 0$,
$$
D_{\bm\lambda}^{-1}\circ \nu\circ\alpha\circ D_{\bm\lambda}  (\mathbf x) $$
uniformly (with respect to $\mathbf x$) converges to $\mathbf x$, and
$$\alpha\circ D_{\bm\lambda}  (\mathbf x)$$  to $\mathbf 0$. 
Thus 
$$\phi (D_{\bm\lambda}^{-1}\circ \nu\circ\alpha\circ D_{\bm\lambda}  (\mathbf x)) 
g(\alpha\circ D_{\bm\lambda}  (\mathbf x))$$ uniformly converges to
$$\phi (\mathbf x) g(\mathbf 0).$$

Then we have
\begin{align*}
A_{\bm\lambda} =  & det(\alpha)\int_{\mathbf x\in \mathbb R^r} 
 L (\alpha\circ D_{\bm\lambda}  (\mathbf x))\cdot   
\biggl( 
\phi (D_{\bm\lambda}^{-1}\circ \nu\circ\alpha\circ D_{\bm\lambda}  (\mathbf x))
g(\alpha\circ D_{\bm\lambda}  (\mathbf x))-\phi ( \mathbf x) g(\mathbf 0)\biggr) d\mu_{x}\\ &
+  det(\alpha) \int_{\mathbf x\in \mathbb R^r}L (\alpha\circ D_{\bm\lambda}  (\mathbf x))
\phi ( \mathbf x) g(\mathbf 0)d\mu_{x}
\end{align*}
Since the function $$ L (\alpha\circ D_{\bm\lambda}  (\mathbf x))$$ is bounded, we conclude 
that
$$
\displaystyle{\lim_{\bm\lambda \Rsh  \mathbf 0}}A_{\bm\lambda}
 =\displaystyle{\lim_{\bm\lambda \Rsh \mathbf 0}}  \int_{\mathbf x\in \mathbb R^r}
 L (\alpha\circ D_{\bm\lambda}  (\mathbf x)) \phi ( \mathbf x)
 d\mu_{x}$$
Since $\alpha\circ D_{\bm\lambda}$ is  the testing map with the invertible linear map
 $\alpha\circ B$,  $$
\displaystyle{\lim_{\bm\lambda \Rsh  \mathbf 0}}  \int_{\mathbf x\in \mathbb R^r}
 L (\alpha\circ D_{\bm\lambda}  (\mathbf x)) \phi ( \mathbf x)
 d\mu_{x}$$
converges. 
This completes the proof.
\par

\end{proof}

\bigskip

\begin{definition} ( of notations) 
We collect notations for some terms in [3].
\par

(1)  The functional of a distribution $\mathcal F$  is denoted by
\begin{equation}
\int_{\mathcal F}{(\bullet)}.
\end{equation}
\par\hspace{1CC} The notation is extended to the functional of a signed measure on  \par\hspace{1CC} 
characteristic functions of measurable sets, or simply the measurable  \par\hspace{1CC} sets.  It should be stressed that
the notation does not
 involve the volume \par\hspace{1CC} element which is excluded  for distribution.\par
(2) Let $\mathbb R^m$ be equipped 
with coordinates
$x=\{x_1, \cdots, x_m\}$.  Let $V_I$ be an \par\hspace{1CC} $r$ dimensional coordinates plane with multi-index $I$ of length $r$, 
and $$\pi_I: \mathbb R^m \to V_I$$ \par\hspace{1CC} be the projection. Let 
 $V_{I^\diamond}$ be
the  perpendicular  coordinates plane of \par\hspace{1CC} dimension $m-r$ satisfying  $\{I I^\diamond\}=\{1, 2, \cdots, m\}$ 
with concordant \par\hspace{1CC} orientations. 
Let $d\mu^I, d\mu^{I^\diamond}$ be  their Euclidean volume forms
$$dx_{i_1}\wedge\cdots \wedge dx_{i_r}, dx_{i_1^\diamond}\wedge\cdots \wedge dx_{i_{m-r}^\diamond}$$
\par\hspace{1CC} with the matching  orders of the $\wedge$  products.  Throughout this paper  \par\hspace{1CC} 
Euclidean volume forms are used in two
different ways interchangeably:  \par\hspace{1CC} a) as a $C^\infty$ differential form with concordant wedge product, b) as the  \par\hspace{1CC} Lebesgue measure on the Euclidean space.   
 \par

(3)   Let $T$ be a current of degree $r$ with a compact support in $\mathbb R^m$. 
In  [3]  \par\hspace{1CC} (Chapter III, \S 8, p36) $T$   is written as
\begin{equation}
T=\sum_{I} \mathcal T_I d\mu^{I^\diamond}.
\end{equation}
\par\hspace{1CC} We call $\mathcal T_I$ for each $I$ the De Rham distribution of $T$.\par

(4) Continuing from (3), 
  let $\mathcal T_{I_1}$ be one of De Rham distributions among \par\hspace{1CC} finitely many $ \mathcal T_I$.  Then
$(\pi_I)_\ast ( \mathcal T_{I_1} d\mu^{I^\diamond})$ is a current of maximal 
degree \par\hspace{1CC} in the plane $V_{I}$ (where  $I_1,  I$ may not be the same). 
 Hence it is regarded \par\hspace{1CC}  as a distribution in $V_I$ 
 (footnote 2 at p34, [3]), denoted by $$(\pi_I)_\star ( \mathcal T_{I_1})$$ \par\hspace{1CC} and 
called the projection of the De Rham distribution  to $V_I$.   \par\hspace{1CC}  
The projection (with the $\star$ subscript) has an expression, 

\begin{equation}
\int_{(\pi_I)_\star ( \mathcal T_{I_1})} f =\int_{T} \pi_I^\ast (f) d\mu^{I_1}.
\end{equation}
\par\hspace{1CC} for a test function $f\in \mathscr D(V_I)$.

\end{definition}

\bigskip

{\bf Remark} \par
  
Part (1): This  unconventional notation will be used for clarity.\par
Part (2): These notations will be used for simplicity.\par
Part (3): The name will be used for convenience.\par
Part (4): This is the usual projection of currents with omission of volume \par\hspace{3.5 CC} forms.  
\smallskip

By above notations,  the  Radon-Nikodym number can also be expressed in distribution as
\begin{equation} \displaystyle{\lim_{ \lambda\Rsh\mathbf 0}} \int_{ (D_{\bm\lambda}^{-1})_{\ast}(L)}\phi, \quad  or\ \quad  
\displaystyle{\lim_{ \lambda\Rsh\mathbf 0}} {\int_ { L} (D_{\bm\lambda}^{-1})^\ast (\phi)\over det( B )\prod_{i=1}^r\lambda_i}.
\end{equation}
\bigskip

\begin{definition}(Lebesgue current).\quad\par

Let $\mathcal X$ be a differentiable manifold of dimension $m$.  Let $U$, a neighborhood,  $x_1, \cdots, x_m$  coordinates for  $U$ be a chart  in the differential  structure of $\mathcal X$. 
Let
$d\mu^I$ be the Euclidean volume form
\begin{equation}
dx_{i_1}\wedge \cdots\wedge dx_{i_r}\end{equation} of an  $r$-dimensional coordinates  
plane $V_I$ with multi-index $I=(i_1\cdots i_r)$,  $\pi_I: U\to V_I\simeq \mathbb R^{r}$  the projection given by the chart.
 Then  a homogeneous current $T$ of dimension $p$ is called Lebesgue if for each chart $(U, x_1, \cdots, x_m)$ in an atlas and 
each set $\mathcal S$ of forms $\xi\in \mathscr D(U)$ with the same compact support and  locally bounded to order $0$ (p38, [3]),  the following conditions are
 satisfied. \medskip

(a) ( Lebesgue condition)  \quad

The projection $(\pi_I)_\star(\mathcal T_{I_1})$ of each De Rham distribution $\mathcal T_{I_1}$ of $T\wedge \xi$ to each   coordinates plane $V_I$  is 
a signed measure absolutely continuous with respect to the Lebesgue measure (under the chart). 
Furthermore the Radon-Nikodym derivative  ${d(\pi_I)_\star(\mathcal T_{I_1})\over d\mu^I}$(see section 32, [1])  
is  an $L^1$   function with a compact support and bounded as free variables vary in  $V_I$ and  $\mathcal S$.   This is equivalent to the existence of 
 a bounded in both $V_I$ and $\mathcal S$,  compactly supported, Lebesgue integrable  function $\mathcal L_I$ on $V_I$ satisfying 
\begin{equation}
\int_{(\pi_I)_\star(\mathcal T_{I_1})} f=\int_{V_I} \mathcal L_I f d\mu^I\end{equation}
for any test function $f\in \mathscr D(V_I)$.  
The $L^1$ function $\mathcal L_I={d(\pi_I)_\star(\mathcal T_{I_1})\over d\mu^I}$  will be called the Lebesgue function of $T$ or $T\wedge \xi$.
 The formula (2.17) can be combined with (2.14) to have a more direct expression of $\mathcal L_I$,
\begin{equation}
\int_{V_I} \mathcal L_I f d\mu^I=\int_{T\wedge \xi }(\pi_I)^\ast (f) d\mu^{I_1}
\end{equation}
where the index $I_1$ is the index  associated to the   De Rham distribution $\mathcal T_{I_1}$, i.e.
$$T\wedge \xi=\mathcal T_{I_1} d\mu^{I_1^\diamond}+\cdots.$$

 \bigskip

(b) (Radon-Nikodym condition)\quad  

All Lebesgue functions $\mathcal L_I$ of $T$  are of Radon-Nikodym.

\end{definition}

\bigskip

{\bf Remark} There are infinitely many Lebesgue functions of $T$. For instance they all depend on $\xi$ which is 
not reflected in the notation $\mathcal L_I$. 
Definition 2.4   addresses the  values of the density function of two measures at  points.\footnote {
   A Radon-Nikodym derivative evaluated at an a.e. point  is the infinitesimal ratio of two measures, called the density.  }
In integral theory  it is more precisely described as follows. 
\bigskip

\begin{proposition} Assume all notations from Definition 2.4.  
 Then Radon-Nikodym condition  of Definition 2.4 holds if and only if 
\begin{equation}\begin{array}{c}
\displaystyle{\lim_{\bm\lambda\Rsh \mathbf 0}} {1\over det ( B)\prod_{i=1}^r\lambda_i}  \int_{\mathbf v\in V_I} 
\mathcal L_I(\mathbf v) \phi (D_{\bm\lambda}^{-1}(\mathbf v)) d\mu^I\end{array}
\end{equation}
\noindent exists for each  test function $\phi$ and index $I$.  
Furthermore if the Lebesgue function $\mathcal L_I$ is continuous at $\mathbf u$,  
\begin{equation}
RN_{\phi, \mathcal L_I}=\mathcal L_I(\mathbf u)\int_{\mathbf v\in V_I}\phi(\mathbf v)d\mu^I.
\end{equation}

\end{proposition}
\bigskip

\begin{proof} 
After the change of variables
\begin{equation}
D_{\bm\lambda}(\mathbf v) \Rightarrow \mathbf v.
\end{equation}
( replacement of  $D_{\bm\lambda}(\mathbf v) $ with $\mathbf v$ ) the integral (2.4) is the same as (2.19), which
says in distributions
\begin{equation}
\int_{(D_{\bm\lambda}^{-1})_\ast (\mathcal L_I)} \phi=
 {1\over det (B)\prod _{k=1}^r \lambda_k} \int_{\mathcal L_I} {(D_{\bm\lambda}^{-1})^\ast (\phi)}.
\end{equation}
If $\mathcal L_I$ is continuous, since $\mathbf v$ is bounded, we have 
\begin{equation}\begin{array}{c}
\displaystyle{\lim_{\bm\lambda\Rsh \mathbf 0}} \int_{\mathbf v\in V_I} 
\mathcal L_I ( D_{\bm \lambda}(\mathbf v) )\phi(\mathbf v) d\mu^I
= \int_{\mathbf v\in V_I} \displaystyle{\lim_{\bm\lambda\Rsh \mathbf 0}}
\mathcal L_I ( D_{\bm \lambda }(\mathbf v) )\phi(\mathbf v) d\mu^I
\end{array}\end{equation}

Therefore the limit exists and is equal to 
$$
 \int_{\mathbf v\in V_I}  \mathcal L_I(\mathbf u)\phi(\mathbf v) )d\mu^I= \mathcal L_I(\mathbf u)\int_{V_I}\phi(\mathbf v)d\mu^I.$$
 Thus the Radon-Nikodym condition is satisfied. 

\end{proof}

\bigskip

{\bf Remark} However the convergence of (1.1) only concerns the limiting patterns  of densities at the non-continuous points.

\bigskip

   Definition 2.4  is stated in one atlas.  Let's show it is independent of the atlas.

 \bigskip

\begin{proposition}
Definition 2.4 defines an invariant of the $C^\infty$ differential structure.

\end{proposition}

\bigskip

\begin{proof}  We need to prove that  the conditions  (a), (b) of Definition 2.4  are independent of charts.  
Let $T$ be a current of dimension $p$, and $\xi\in \mathscr D(U)$ a  form in a neighborhood
$U$.  Let $U, x=\{x_1, \cdots, x_m\}$ be  a chart called $x$-chart satisfying the conditions of Definition 2.4 for $T\wedge \xi$. 
  Let  $ U, y=\{y_1, \cdots, y_m\}$ be another chart called $y$-chart. 
Let $\nu$ be the transition map from $x$-chart to $y$-chart. 
Let $V_I$ be an $r$ dimensional $x$ coordinates plane,  $V_J$ be an $r$ dimensional $y$ coordinates plane,
$d\mu_x^I, d\mu_y^J$ be the volume forms of  the coordinates planes $V_I, V_J$ respectively. Let
\begin{equation}
d\mu_y^{J^\diamond}=\sum_{K} g_{JK}(\mathbf x)d\mu_x^{K^\diamond}
\end{equation}
where $g_{JK}$ is the entry of the Jacobian matrix $J_{y\to x}$ from $y$-chart to $x$-chart, and $K$ is a multi-index of length $r$. 
Let $\pi_J: U\to V_J$ be the projection through $y$-chart, and  $\pi_I: U\to V_I$ the projection through the $x$-chart. 
 We may assume the projection map (through $y$-chart)
$$\begin{array}{ccc}
\nu_{IJ}:  V_I &\rightarrow & V_J \end{array}$$ is a diffeomorphism that preserves the orientation.   
Now we  fix $J$ index of length $r$.
On $U$, we have the sum
\begin{equation}
T\wedge \xi= \sum_K F_K(\mathbf y) d\mu_y^{K^\diamond}\end{equation}
where $ F_K(\mathbf y)$ is a De Rham distribution on $U$, and $K$ is the multi-index of length $r$. 
Then $supp(F_K(\mathbf y) )$ is bounded,  since 
$T\wedge \xi$ has a compact support.   Then for two fixed indexes $J, K$ of length $r$, 
 $$F_K(\mathbf y) d\mu_y^{J^\diamond}$$
is a current on $U$ of dimension $r$. Through $x$-chart it has the following decomposition

$$ F_K(\mathbf y) d\mu_y^{J^\diamond}=\sum_{I} \mathcal D_I$$
where \begin{equation}
\mathcal D_I= F_K(\nu (\mathbf x)) g_{JI}(\mathbf x) d\mu_x^{I^\diamond}\end{equation}
is a current of dimension $r$ on $U$, 
and $I^{\diamond}$ is a multi-index of length $m-r$.  Note: $F_K(\nu (\mathbf x))$ is the distribution 
$$(\nu^{-1})_\ast ( F_K(\mathbf y)).$$
This notation for push-forwards of distributions will be used alternately with the conventional notations throughout, 
but this is referred to as the change of variables. 
\bigskip

There is a commutative diagram
\begin{equation}
\begin{tikzcd}
 &  U \arrow[dl, phantom , shift left]
\arrow[dl, "\pi_I" above left]  \arrow{dr}{\pi_J} \\
V_I  \arrow{rr}{\nu_{IJ}} && V_J.
\end{tikzcd}
\end{equation}
Then we have 
\begin{align*}
 &(\pi_J)_\star ( F_K(\mathbf y)) \\ &=
  \sum_I  (\pi_J)_\ast    ( \mathcal D_I ) \\
&=\sum_I(\nu_{IJ})_\ast\circ (\pi_I)_\ast (\mathcal D_I)
\end{align*}
(Note: $\star$, $\ast$ are two different operators.). 
Therefore for distributions we have 
\begin{equation}
 (\pi_J)_\star ( F_K(\mathbf y))=
 \sum_I(\nu_{IJ})_\ast \circ (\pi_I)_\ast (\mathcal D_I).
\end{equation}
 
   Let's calculate $\mathcal D_I$.
Let
\begin{align*}
 &T\wedge \xi=\sum_P G_P(\mathbf x) d\mu^{P^\diamond}_x\\
 &(\text {Note:  $G_P(\mathbf x)$ is a distribution})\\
&=\sum_K \sum_P G_P (\nu^{-1} (\mathbf y)) g^{-1}_{PK}(\mathbf y) d\mu_y^{K^\diamond},\\
\end{align*}
where $ g^{-1}_{PK}$ stands for the entry of the Jacobian matrix, $J_{x\to y}$.
Hence 
$$ F_K(\mathbf y)=\sum_P G_P (\nu^{-1} (\mathbf y)) g^{-1}_{PK}(\mathbf y). $$
Then 
 \begin{equation}
\mathcal D_I= \sum_P G_P (\mathbf x) g^{-1}_{PK}( \nu (\mathbf x)) g_{JI}(\mathbf x) d\mu_x^{I^\diamond}.\end{equation}

Since $T$ is Lebesgue in $x$-chart, it satisfies both  conditions of Definition 2.4 in $x$-chart, therefore
$(\pi_I)_\ast (\mathcal D_I)$
is a bounded, compactly supported $L^1$ function of Radon-Nikodym on $V_I$ (in $x$-chart). 
Due to Proposition 2.2, so is  
$$
(\nu_{IJ})_\star \circ (\pi_I)_\ast (\mathcal D_I)$$  on $V_J$.  Hence its sum over finitely many $I$, 
 $$\sum_I(\nu_{IJ})_\ast\circ (\pi_I)_\ast (\mathcal D_I)$$ is also  a bounded, compactly supported $L^1$ function of Radon-Nikodym on $V_J$ (which is in $y$-chart).  By the formula (2.28) we complete the proof.

\end{proof}

\bigskip

\subsection{Examples}

 We work with singular chains of regular cells. 
\bigskip

\begin{proposition}
 Let $c$ be a cell. Then $c$ is Lebesgue. Furthermore chains of  cells are Lebesgue. So singular cycles are Lebesgue.  \end{proposition}

\bigskip

\begin{proof}  It suffices to work in one chart.  So we  assume $\mathcal X=U=\mathbb R^m$ is equipped 
with the standard chart $\mathscr A$ (a basis for the linear space) with coordinates $(x_1, \cdots, x_m)$. 
We may assume  the cell $c$ is represented by 
a diffeomorphism extended to a neighborhood $\mathcal K$
\begin{equation}\begin{array}{ccc}
h: \mathcal K &\rightarrow & h(\mathcal K)\subset U\\
\quad \cup & & \cup\\
\quad \Delta &\rightarrow & h(\Delta)
\end{array}\end{equation}
where $\Delta$ is an $r$-dimensional  polyhedron equipped with the Lebesgue measure  $d\mu$.
Let $\xi$ be a test form in $\mathscr D(U)$. 
Let $V_I\simeq \mathbb R^r$ be an $r$-dimensional  coordinates plane. We denote  projection  $U\to V_I$ by $\pi_I$. Let $d\mu^J$ be the 
Euclidean volume form of another $r$-dimensional coordinates plane. Then by the formula (2.18) the projection of   
a De Rham distribution of $c\wedge \xi$ to $V_I$ is the functional, 
\begin{equation}\begin{array}{ccc}
\mathcal F: \phi &\rightarrow & \int_{c\wedge \xi} \phi(\mathbf x) d\mu^J
\end{array}\end{equation}
where $\phi(\mathbf x)=\pi_I^\ast ( \phi(\mathbf v))$  for $\phi(\mathbf v)\in \mathscr D(V_I)$ ( Note: the index $J$ is associated with the De Rham distribution).

Notice that
\begin{equation}\begin{array}{c}
| \int_{c\wedge \xi} \phi(\mathbf x) d\mu^J|\\
=|\int_c \xi\wedge \phi(\mathbf x) d\mu^J | \\
=|\int_{\Delta^p} h^\ast (\xi\wedge \phi(\mathbf x) d\mu^J)|\\
\leq C ||\phi||_{0, K}
\end{array}\end{equation}
 where $C$ is a constant and $||\phi||_{0, K}$ is semi-norm with the compact support $K=supp(\phi(\mathbf v))$. 
Since the inequality (2.32) holds for all compact set $K$ supporting the $\phi$, by Proposition 2.1.8, 1.3.11, [4], $\mathcal F$ is a measure. 
So if we let $\phi$ be a characteristic function $\chi(E)$ of a subset set $E\subset V_I$ of Lebesgue measure $0$, the inequality (2.32) becomes
\begin{equation}\begin{array}{c}
\int_{|\mathcal F|} \chi(E) \leq C \int _  {h^{-1}(E) }d\mu =0.
\end{array}\end{equation}
Hence $\mathcal F$ is absolutely continuous with respect to the Lebesgue measure of $V_I$. 
Next we estimate the Radon-Nikodym derivative  a.e.
\begin{equation}
\displaystyle{\lim_{\epsilon\to 0}}{|\int_{\Delta^p} h^\ast (\xi\wedge \chi(B_\epsilon) d\mu^J)|\over d\mu_I|_{\chi(B_\epsilon)}}
\end{equation}
where $B_\epsilon$ is a bounded domain  in $V_I$ of radius $\epsilon$.  Notice that 
$$
\biggl|\int_{\Delta^p} h^\ast (\xi\wedge \chi(B_\epsilon) d\mu^J)\biggr|\leq C||\xi||_{0, K} \int_{B_\epsilon} d\mu_I.
$$
Therefore \begin{equation}
\biggl|{\int_{\Delta^p} h^\ast (\xi\wedge \chi(B_\epsilon) d\mu^J)\over \int_{B_\epsilon} d\mu_I }\biggr|\leq 
C ||\xi||_{0, K}.
\end{equation}
Hence the Lebesgue function ${d\mathcal F\over d\mu^I}$ is bounded when $\xi$ is locally bounded to order $0$. 
This shows $c$ satisfies the Lebesgue condition. 

\bigskip

 Let's now prove the Radon-Nikodym condition.  We may assume $\mathbf u=\mathbf 0$ and $B=identity$. 
Let $\phi(\mathbf v)$ be the test function on $V_I$. 
Recall the projection $\pi_I: U\to V_I$, $D_{\bm\lambda}$ the testing map which is a block-wise scalar multiplication (see Definition 2.1).  
Since being Lebesgue current is independent choice of coordinates, we may choose a coordinates system so that  the composition 
$$P: \mathcal K\to U\to V_I$$ is a diffieomorphism.  Then by using the formula (2.18), the Radon-Nikodym number is the limit 
\begin{equation}
 \displaystyle{\lim_{\bm\lambda \Rsh 0}}
\int_{c\wedge \xi }   \pi_I^\ast \biggl (  {\phi( D_{\bm\lambda}^{-1} (\mathbf v))\over det (\mathbb D_{\bm\lambda} B)}
\biggr) d\mu^J,
\end{equation} 
 where $ J$ is an arbitrary multi-index of length $r$, and $D_{\bm\lambda}(\bullet)$ is the affine transformation as in (2.3)
( the integrand is the $C^\infty$ form on $\mathbb R^r$). Because the integrand can absorb the $C^\infty$ form $\xi$,   
we may  assume $\xi$ has degree $0$ and has value $1$ on $\bar c$, and $I=J$.     Then after the change of 
variables 
$$D_{\bm\lambda}^{-1} (\mathbf v) \Rightarrow \mathbf v, $$

the integral (2.36) is the evaluation of distributions on the plane $V_I\simeq \mathbb R^r$, 
\begin{equation} 
\int_{  D_{\bm\lambda}^{-1}( P(\Delta)) } \phi(\mathbf v), 
\end{equation}
where $ D_{\bm\lambda}^{-1} ( P(\Delta) ) $  is a cell for each $\bm\lambda$.

\bigskip

Next  we use measure theory. 
\bigskip

\begin{ass}
The limit of distributions
$$ D_{\bm\lambda}^{-1} ((P (\Delta))) $$
as $|\bm\lambda|\Rsh 0$ exists.

\end{ass}

\begin{proof} of claim 2.8:  First let's have a definition in analysis.\bigskip

Let $\mathbb R^{k_1}$ be a subspace of $\mathbb R^r$ with the natural number $k_1$.  So
with any coordinate's chart ( linear basis), there is a direct sum decomposition  \begin{equation}
\mathbb R^r= \mathbb R^{k_1}\oplus \mathbb R^{k_2}.
\end{equation}
   Let $ W\subset \mathbb R^r$ be a bounded measurable set. 
Let $W_{q_2}$ be the
the collection of points $x$ in $\mathbb R^{k_1}$ such that  $(x, q_2)$  lies in $W$.  Let 
$S^{k_1-1}$ be the unit sphere of $\mathbb R^{k_1}$.  
\bigskip

\begin{definition}   Let $B_2$ be a small ball of $\mathbb R^{k_2}$ around the origin. 
We say $W$ is a growing set along $\mathbb R^{k_1}$
 if  there exist $B_2$ and $\epsilon>0$  such that for  a.e point $(v_1, q_2)\in S^{k_1-1}\times B_2$
 the line segment 
\begin{equation} L_\epsilon=\{ tv_1:  \ 0<t<\epsilon\}\subset \mathbb R^{k_1}\end{equation}
either lies in $W_{q_2}$ or does not meet $W_{q_2}$, where the a.e for the unit vector $v_1$ is the abbreviation 
of  ``almost everywhere"  in the spherical measure on the unit sphere $S^{k_1-1}\subset \mathbb R^{k_1}$. \end{definition}\bigskip

We divide the proof into two steps.\medskip

\noindent {\it Step 1}:  We assume the  decomposition (2.38) satisfies   \begin{assumption}
 $P(a)=(\mathbf 0, \mathbf 0)$, \begin{equation}
T_{(\mathbf 0, \mathbf 0)} (S^{k_1-1}\times \{\mathbf 0\})\not\subset  P_\ast (T_{a}\partial \bar\Delta^r)
\end{equation}
where $T_{a}\partial \bar\Delta^r$ is the union of finitely many tangent planes to each face 
 in $T_{a}\mathbb R^r$.

\end{assumption}
\bigskip

 Due to Assumption 2.10, for an a.e. $(v_1, q_2)\in S^{k_1-1}\times B_2$, 
\begin{equation}
 v_1\not\in P_\ast (T_{(\mathbf 0, \mathbf 0)}\partial \bar\Delta^r).\end{equation}
Thus $P(\bar\Delta^r)$ is a growing set which is divided into two pieces:
$$P(\bar\Delta^r)=\mathcal E_{1} \cup \mathcal E_{2}$$
where 
$\mathcal E_{1}$ lies in the union of line segments $\cup (L_{\epsilon}\times \{q_2\})$
and $\mathcal E_{2}$ is the complement of $\mathcal E_1$.
Let $D_{\lambda, k_1}$ be the linear transformation of the direct sum of 
the scalar multiplication map by $\lambda>0$ on $\mathbb R^{k_1}$ and identity on
$\mathbb R^{k_2}$. Because the formula (2.39) is satisfied,
\begin{equation}
\displaystyle{\lim_{\lambda\to 0}} \int_{D_{\lambda, k_1}(\mathcal E_{2})}\phi(\mathbf v)=0.
\end{equation}

For the set $ D_{\lambda, k_1} (\mathcal E_{1})$, we defined a sequence of measurable sets 
\begin{equation}
A_n=D_{{1 \over n}, k_1} (\mathcal E_{1}).
\end{equation}
Due to the condition (2.39), we have the sequence of cells
\begin{equation}
A_1\subset A_2 \subset \cdots\subset A_n\subset  \cdots\end{equation}
Let $A_\infty=\cup_{n=1}^\infty A_n$. 
 Hence the measurable sets $A_n$ in measure converges to $A_\infty$. 
So $D_{\lambda, k_1}(P(\bar\Delta^r))$ converges to $A_\infty$ as $\lambda\to 0$.

\bigskip

Recall the zigzag path is a division of basis $\mathbf e_{1}, \cdots, \mathbf e_{r}$ into $l$ groups 
 with an order:  $j_1$ group , $j_2$ group, $\cdots$,  $j_l$ group.  
We define the subspaces of $\mathbb R^r$ spanned by the basis of $j_i$ groups by 
$$\mathbb R^{r_j}, j=1, \cdots, l. $$
Then
\begin{equation}
\mathbb R^r=\mathbb R^{r_1}\oplus \cdots\oplus \mathbb R^{r_l}.
\end{equation}
Now we inductively repeat the step 1 for each subspace $\mathbb R^{r_j}$ to obtain
the growing set along $\mathbb R^{r_j}$. This leads to the convergence:  
\begin{equation}
\displaystyle{\lim_{|\bm\lambda|\Rsh 0}}\int_{ D_{\bm\lambda}^{-1}(P(\Delta) ) }\phi(\mathbf v)
\end{equation} exists. 
More specifically the distributions converge to a finite, measurable set $E$, i.e.  
$$\displaystyle{\lim_{|\bm\lambda|\Rsh 0}} {D_{\bm\lambda}^{-1}(P(\Delta) ) }=E.$$

\bigskip

\noindent 
{\it Step 2}:  In this step, we remove Assumption 2.10 which amounts
to have a linear transformation of the original coordinates of $\mathbb R^r$. So let 
$\mathbb R^r$ be equipped with  standard basis $\mathscr A$.  
We choose a family of orthonormal linear transformation $G_{t}\in SU(r), t\in\mathbb R$ such that 
$$\displaystyle{\lim_{t\to 0}}G_t=identity.$$
So if $G_t: \mathbb R^{r}\to \mathbb R^{r}$ is a linear transformation, $G_t(\mathscr A), $ for $t\neq 0$ satisfies the assumption 2.10. Then  the limit as $|\bm\lambda|\Rsh 0$
\begin{equation}
 \int_{D_{\bm\lambda}^{-1}(t)(P(\Delta) ) }\phi( \mathbf v)
\end{equation} exists for each fixed $\bm\lambda$, 
where $ D_{\bm\lambda}^{-1}(t)$ is the  testing map under the coordinates of the chart $G_t(\mathscr A)$.
Notice the original integral (2.37) is 
\begin{equation}
 \int_{D_{\bm\lambda}^{-1}(0)(P(\Delta) ) }\phi( \mathbf v).
\end{equation}

First we notice $\phi$ is $C^\infty$ and has a bounded support. This  leads to an  estimate
\begin{equation} \biggl| \int_{D_{\bm\lambda}^{-1}(t)(P(\Delta) ) }\phi( \mathbf v)-
\int_{D_{\bm\lambda}^{-1}(0)(P(\Delta) ) }\phi( \mathbf v)\biggr|
\leq C' |t|
\end{equation}
for some constant $C'$ that is independent of $\bm\lambda$. 
Since  $$\displaystyle{\lim_{|\bm\lambda|\Rsh 0}}\int_{D_{\bm\lambda}^{-1}(t)(P(\Delta) ) }\phi( \mathbf v)$$
exists for each $t$, so does $$\displaystyle{\lim_{|\bm\lambda|\Rsh 0}}\int_{D_{\bm\lambda}^{-1}(0)(P(\Delta) ) }\phi( \mathbf v)$$

We complete the proof of Claim 2.8.
\end{proof}

\bigskip

By the linearity, the existence is extended to all chains and cycles 
\end{proof}
\bigskip

\begin{proposition} Let $\omega$ be a $C^\infty$ form. Then $\omega$ is Lebesgue. \end{proposition}
\bigskip

\begin{proof} We may prove it locally. So let $\mathcal X=\mathbb R^m$. Let $m-p=deg(\omega)$.
Let $\xi$ be any test form on $\mathbb R^m$ of degree $p-r$.
Let $$x_1, \cdots, x_{r}, x_{r+1},\cdots,  x_{p},  x_{p+1}, \cdots, x_{m}$$
 be a coordinates  chart. Let
$V_I$ have coordinates plane of components $x_1, \cdots, x_r$.   
For the simplicity, we may assume
\begin{equation}
\omega=M(\mathbf x) dx_{p+1}\wedge\cdots\wedge dx_m
\end{equation}
\begin{equation}
\xi= N(\mathbf x)dx_{r+1}\wedge\cdots\wedge dx_p
\end{equation}
We obtain that
the Lebesgue function of $ \omega\wedge \xi$ is
the fibre integral 
\begin{equation} 
\int_{(x_{r+1}, \cdots, x_m)\in \mathbb R^{m-r}}M(\mathbf x)N(\mathbf x)dx_{p+1}\wedge\cdots\wedge dx_m\wedge dx_{r+1}\wedge\cdots\wedge dx_p\end{equation}
which is
a $C^\infty$ function of $x_1, \cdots, x_r$ in the $V_I$ plane.   Since the Lebesgue function is $C^\infty$, by Proposition 2.5, 
the Radon-Nikodym condition is  satisfied 
\end{proof}

\bigskip

\begin{proposition} Let $\mathcal Y$ be another $C^\infty$ manifold. 
If currents $T_1, T_2$ are Lebesgue  in $\mathcal X, \mathcal Y$ respectively, so is 
$T_1 T_2$  in $\mathcal X\times \mathcal Y$, where $T_1T_2$ is the tensor product of currents defined 
as in \S 12, [3].

\end{proposition}

\bigskip

First we need some lemmas for the proposition.
\bigskip

\begin{lemma}
We resume the set-up of Definition 2.4. In particular 
$U$ is a neighborhood of a chart in manifold $\mathcal X$.
If \begin{equation}
\begin{array}{c}
\phi\in \mathscr D(U\times \mathbb R^n),\\
\end{array}\end{equation}
 
(1) then  Radon-Nikodym number $RN_{\phi, \mathcal L_I}(\mathbf y)$,  which is a function of \par\hspace {1CC}
$\mathbf y\in \mathbb R^n$ lies in $ \mathscr D( \mathbb R^n)$. \par

(2) The convergence 
\begin{equation}\begin{array}{c}
\displaystyle{\lim_{\bm \lambda\Rsh  0}} \int_{\mathbf v\in V_I} 
\mathcal L_I ( D_{\bm \lambda }(\mathbf v) )\phi(\mathbf v, \mathbf y) d\mu^I\end{array}\end{equation}
\par\hspace {1CC}  is uniform for  the bounded variable $\mathbf y\in \mathbb R^n$.

\end{lemma}

\bigskip

Before to prove them, we give a definition 
in analysis. 
\bigskip

\begin{definition} ( ``$\epsilon$-$\delta$" for ordered limits)  Let  $\bm\lambda$ have the grouping as in Definition 2.1. 
Let $N_{\bm\lambda}$ be  a real function of $\bm\lambda$. Assume the zigzag limit
 $$\displaystyle {\lim_{\lambda_{j_l}\to 0}}\cdots \displaystyle {\lim_{\lambda_{j_1}\to 0}}N_{\bm\lambda}$$
exists. Then for each $\epsilon>0$, there exist a real number $\delta$ and  the  sequence
$$\lambda_{j_1}, \lambda_{j_2}, \cdots, \lambda_{j_l}$$ satisfying that  $\lambda_{j_k}$ is determined by all $\lambda_{j_{k'}}, k'< k$
with the last bound $|\lambda_{j_l}|\leq \delta$
such that  a sequence of inequalities 
$$|N_{(\lambda_{j_l}, \cdots, \lambda_{j_k}, 0, \cdots, 0)}-N_{(\lambda_{j_l}, \cdots, \lambda_{j_{k+1}}, 0, \cdots, 0 )}|\leq \epsilon$$
hold, where $N_{(\lambda_{j_l}, \cdots, \lambda_{j_i}, 0, \cdots,0 )}$ for all $i$ are defined as the limits:
$$\displaystyle {\lim_{\lambda_{j_{i-1}}\to 0}}\cdots \displaystyle {\lim_{\lambda_{j_1}\to 0}}N_{\bm\lambda}.$$

We  call such $\bm\lambda$  a point on the zigzag path, and denote the condition on $\bm\lambda$ by
$$|||\bm\lambda|||\prec\epsilon .$$
So the usual ``$\epsilon-\delta$" statement in the analysis can be applied to the zigzag limit 
with the notion 
$|||\bm\lambda|||\prec\epsilon$. 
\end{definition}
\bigskip

\begin{proof} 
\quad\par
(1)

Express the  Radon-Nikodym number  as a function of $\mathbf y$,
$$RN_{\phi,  \mathcal L_I}(\mathbf y)$$
where $\mathbf y\in \mathbb R^n$ and $\mathcal L_I$ is a Lebesgue function.   Let $\mathbf e_i, i=1, \cdots, n$ be a basis for $\mathbb R^n$. 
Let $\Delta y$ be a real number, $\mathbf y=\sum_{i=1}^n y_i\mathbf e_i\in \mathbb R^n$.
Let's consider the number
\begin{align*}
A_h  &={ RN_{\phi,  \mathcal L_I}(\mathbf y+\Delta y\mathbf e_i)-RN_{\phi,  \mathcal L_I}(\mathbf y)\over \Delta y}- RN_{{\partial \phi(\mathbf y)\over \partial y_i}, \mathcal L_I}(\mathbf y) \\&
 =\displaystyle{\lim_{|\bm\lambda|\Rsh \mathbf 0}}
\int_{\mathbf x\in V_I} \mathcal L_I(D_{\bm\lambda}(\mathbf x)) \biggl ({\phi\bigl (D_{\bm\lambda}(\mathbf x), \mathbf y+\Delta y\mathbf e_i
\bigr)-
\phi\bigl(D_{\bm\lambda}(\mathbf x), \mathbf y\bigr)\over \Delta y}-{\partial \phi\bigl(D_{\bm\lambda}(\mathbf x),\mathbf y\bigr)\over \partial y_i}\biggr)d\mu^I
\end{align*}
Since $\phi$ is $C^\infty$ with a compact support,  
$${\phi\bigl (D_{\bm\lambda}(\mathbf x), \mathbf y+\Delta y\mathbf e_i
\bigr)-
\phi\bigl(D_{\bm\lambda}(\mathbf x), \mathbf y\bigr)\over \Delta y}-{\partial \phi\bigl(D_{\bm\lambda}(\mathbf x),\mathbf y\bigr)\over \partial y_i}$$
as $h\to 0$ uniformly (with respect to $\bm\lambda, \mathbf x$) converges to $0$.  Together with the
bounded  $\mathcal L_I(D_{\bm\lambda}(\mathbf x))$, we have 
 $$\displaystyle{\lim_{\Delta y\to \mathbf 0}}A_h=0.$$
Hence $RN_{\phi,  \mathcal L_I}(\mathbf y)$ is differentiable and 
\begin{equation}
{\partial RN_{\phi,  \mathcal L_I}(\mathbf y)\over \partial y_i}=RN_{ {\partial \phi(\mathbf y)\over \partial y_i}, \mathcal L_I}(\mathbf y)
\end{equation}
By the iteration, $RN_{\phi,  \mathcal L_I}(\mathbf y)$ is $C^\infty$. Since $\phi(\mathbf x, \mathbf y)$ is both bounded and  compactly supported, 
so is $RN_{\phi,  \mathcal L_I}(\mathbf y)$. 

\bigskip

(2)  Let's continue from part (1). By Theorem 6, Chapter II, \S 7, [3], there is a sequence of test functions
\begin{equation}
\psi^n_1(\mathbf v)\in \mathscr D(V_I), \psi^n_2(\mathbf y)\in \mathscr D(\mathbb R^n)
\end{equation}
 
such that $$\psi^n_1(\mathbf v)  \psi^n_2(\mathbf y)\to \phi(\mathbf v, \mathbf y)\ as\ n\to \infty$$
in $\mathscr D(V_I\times \mathbb R^n)$.
Thus for any $\epsilon>0$, since $\mathcal L_I$ is bounded
there is an $N$ such that
\begin{equation}
|\int_{\mathbf v\in V_I} 
\mathcal L_I ( D_{\bm \lambda }(\mathbf v) )\biggl(\psi^N_1(\mathbf v) \psi^N_2(\mathbf y)
-\phi(\mathbf v, \mathbf y)\biggr) d\mu^I|\leq \epsilon
\end{equation}
for all $\bm\lambda$.  Taking the limit $|\bm\lambda|\Rsh 0$, we have

\begin{equation}
|RN_{\psi_1^N , \mathcal L_I}\psi_2^N(\mathbf y)-RN_{\phi,  \mathcal L_I}(\mathbf y)|\leq \epsilon 
\end{equation}
Now for this fixed $N$,  $\underset{\mathbf y}{Max}   ( |\psi_2^N(\mathbf y)|$ is a fixed number. Hence 
there is a $\delta$ such that for when $\bm\lambda$ is on zigzag path and  $|||\bm\lambda|||\prec\epsilon$, 
\begin{equation}\begin{array}{c}
|RN_{\psi_1^N, \mathcal L_I}\psi_2^N(\mathbf y)-
\int_{\mathbf v\in V_I}\mathcal L_I ( D_{\bm \lambda }(\mathbf v) )\psi^N_1(\mathbf v)\psi_2^N(\mathbf y)
 d\mu^I|\leq \epsilon .\end{array}\end{equation} 

So for $\bm\lambda$ on a zigzag path with $|||\bm\lambda|||\prec\epsilon$, we have four numbers varied with $\mathbf y$ and  differed by $\epsilon$:
\begin{equation}\begin{array}{c}
\int_{\mathbf v\in V_I} 
\mathcal L_I ( D_{\bm \lambda }(\mathbf v) \phi(\mathbf v, \mathbf y) d\mu^I, 
\\
 \downarrow \\
\int_{\mathbf v\in V_I}\mathcal L_I ( D_{\bm \lambda }(\mathbf v) )\psi^N_1(\mathbf v) \psi^N_2(\mathbf y)d\mu^I \\
\downarrow\\ 
RN_{\psi_1^N, \mathcal L_I}\psi_2^N(\mathbf y)\\
\downarrow \\
RN_{\phi, \mathcal L_I}(\mathbf y) 

.\end{array}\end{equation}
Hence the convergence 
\begin{equation}
 \int_{\mathbf v\in V_I} 
\mathcal L_I ( D_{\bm \lambda }(\mathbf v) \phi(\mathbf v, \mathbf y)  d\mu^I\to RN_{\phi,  \mathcal L_I}
\end{equation} as $|\bm\lambda|\Rsh\mathbf 0$
is uniform.

\end{proof}

\bigskip

\bigskip

 \begin{proof} of Proposition 2.10:   It suffices to work with $\mathcal X\simeq \mathbb R^m$ 
and $ \mathcal Y\simeq \mathbb R^n$. By Proposition 2.6, it suffices to work with one chart. So
we assume $\mathcal X$ is covered by a single chart $x_1, \cdots, x_m$ and $\mathcal Y$ is covered by a single chart $y_1, \cdots, y_n$. 
Then  $x_1, \cdots, x_m, y_1, \cdots, y_n$ give  a  chart  to $\mathcal X\times \mathcal Y$, called
$x$-$y$ chart.   Let
\begin{equation}
\mathbf e_1, \cdots, \mathbf e_m, \ and \ \mathfrak b_1, \cdots, \mathfrak b_n
\end{equation}
be the  bases for $\mathcal X, \mathcal Y$. First we describe the coordinates planes to which the currents will be projected. 
 For the clarity, we'll use the indexes 
in the following convention.\bigskip

(I) Single index denotes objects from each individual manifold $\mathcal X$ or $\mathcal Y$.\par
\hspace{1 CC}
Indexes $p, k$ denote the objects in $\mathcal X$, $q, l$ in $\mathcal Y$.\par
(II)  Double indexes denote the objects from $\mathcal X\times \mathcal Y$.\par

(III) $d\mu_{\bullet}, d\mu_{\bullet, \bullet}$ are the Lebesgue measures under the charts.
\bigskip

 Using bases, we let

\begin{align*}
&V_k=span(\mathbf e_{1}, \cdots, \mathbf e_{k}),\\
&V_{p-k}=span(\mathbf e_{k+1}, \cdots, \mathbf e_{p}),\\ 
&V_{p}=span(\mathbf e_1, \cdots,  \mathbf e_{p})\\
&V_l=span(\mathfrak b_{1}, \cdots, \mathfrak b_{l}),\\
&V_{q-l}=span(\mathfrak b_{l+1}, \cdots, \mathfrak b_{q}),\\
&V_q=span (\mathfrak b_1, \cdots, \mathfrak b_q). \\
& V_{k, l}=V_k\times V_l\\
 & ...
\end{align*}

\bigskip

 Recall $T_1, T_2$ are currents. Let's assume $dim(T_1)=p, dim(T_2)=q$.  We may assume the  form $\xi$ is in the format  
\begin{equation}
\xi=\zeta(\mathbf x, \mathbf y) d\mu_{p-k, q-l}
\end{equation} with the function 
$\zeta\in \mathscr D(\mathcal X\times \mathcal Y)$.
Let $\xi_x\in \mathscr D(\mathcal X), \xi_y\in \mathscr D(\mathcal Y)$ be functions such that they are equal to 1 on the projections of
$supp(\zeta)$ to $\mathcal X, \mathcal Y$.  
We denote $\xi_x T_1$ by $T_x$ and $\xi_y T_2$ by $T_y$. They all have compact supports.
Then \begin{equation}
(T_1T_2)\wedge \xi=(T_x\wedge d\mu_{p-k}) ( T_y\wedge d\mu_{q-l}) \zeta(\mathbf x, \mathbf y) .\end{equation}
Let $d\mu_x, d\mu_y$ be the Euclidean volume forms of arbitrary coordinates planes of dimensions
$p-k$,  $q-l$ in $x$ chart and $y$ chart respectively. Then by the formula (2.18), the projection $\mathcal T_{k, l}(\xi)$ of a De Rham distribution
of current $(T_1 T_2)\wedge \xi$ to
$V_{k, l}$ is the functional
\begin{equation}\begin{array}{ccc}
\mathcal T_{k, l}(\xi): \phi &\rightarrow &  \int_{(T_x\wedge d\mu_{p-k})( T_y\wedge d\mu_{q-l}) \zeta(\mathbf x, \mathbf y)}
\phi(\mathbf x_k, \mathbf y_l) d\mu_x\wedge d\mu_y\\
 & &= \int_{T_y\wedge d\mu_{q-l}}\biggl( \int_{(T_x\wedge d\mu_{p-k})\zeta(\mathbf x, \mathbf y)} \phi(\mathbf x_k, \mathbf y_l) d\mu_x\biggr) d\mu_y.
\end{array}
\end{equation}
where $\phi(\mathbf x_k, \mathbf y_l)$ is a test function on the coordinates plane $V_{k, l}$. 
By Theorem 6, Chapter II, \S 7, [3], there are sequences of test functions
\begin{equation}
\zeta_x^n (\mathbf x), \zeta_y^n(\mathbf y)
\end{equation}
on $\mathcal X, \mathcal Y$ respectively 
such that $$ \zeta_x^n (\mathbf x) \zeta_y^n(\mathbf y)\to \zeta, \ as\ n\to \infty$$
Then  for any natural number $n$, we have
\begin{equation}\begin{array}{c}
\int_{\mathcal T_{k, l}(\xi)}\phi=\int_{T_y\wedge d\mu_{q-l}}\biggl( \int_{(T_x\wedge d\mu_{p-k})\bigl(\zeta(\mathbf x, \mathbf y)-
\zeta_x^n (\mathbf x)\zeta_y^n(\mathbf y)\bigr)
} \phi(\mathbf x_k, \mathbf y_l) d\mu_x\biggr) d\mu_y\\
+\int_{T_y\wedge d\mu_{q-l}\zeta_y^n(\mathbf y)}\biggl( \int_{T_x\wedge d\mu_{p-k}
\zeta_x^n (\mathbf x)
} \phi(\mathbf x_k, \mathbf y_l) d\mu_x\biggr) d\mu_y\end{array}
\end{equation}
Now we let $\mathcal L_k^n(\mathbf x_k), \mathcal L_l^n (\mathbf y_l)$ be the Lebesgue functions of
De Rham distributions of the currents
$$T_x\wedge d\mu_{p-k}\zeta_x^n (\mathbf x), \quad\quad T_y\wedge d\mu_{q-l}\zeta_y^n(\mathbf y)$$
on $V_k, V_l$ respectively.  The the Lebesgue condition implies 
they are bounded for all $n$, i.e. there is constant $M$ such that
\begin{equation}\begin{array}{c}
|\mathcal L_k^n(\mathbf x_k)|\leq M\\
|\mathcal L_l^n (\mathbf y_l)|\leq M.
\end{array}\end{equation}
By part (1), Lemma 2.13, 
\begin{equation}\begin{array}{c}
\int_{T_y\wedge d\mu_{q-l}\zeta_y^n(\mathbf y)}\biggl( \int_{T_x\wedge d\mu_{p-k}
\zeta_x^n (\mathbf x)
} \phi(\mathbf x_k, \mathbf y_l) d\mu_x\biggr) d\mu_y\\
=\int_{V_{k, l}}\mathcal L_k^n(\mathbf x_k)\mathcal L_l^n (\mathbf y_l)\phi(\mathbf x_k, \mathbf y_l) d\mu_k\wedge d\mu_l
.\end{array}\end{equation}
On the other hand, there is a sequence of numbers 
$a_n\to +\infty$ as $n\to \infty$ such that the set of forms $$a_n ((\zeta(\mathbf x, \mathbf y)-
\zeta_x^n (\mathbf x)\zeta_y^n(\mathbf y))$$for all $n$
is locally bounded to order $0$.  Hence two Lebesgue functions of two currents
$$T_y\wedge d\mu_{q-l}, \quad a_n(T_x\wedge d\mu_{p-k})\bigl(\zeta(\mathbf x, \mathbf y)-
\zeta_x^n (\mathbf x)\zeta_y^n(\mathbf y)\bigr)$$
on $V_l, V_k$ are bounded (the  Lebesgue function on $V_k$ is dependent of $\mathbf y$, but it also bounded for all $\mathbf y$.). 
So the sequence of numbers  $$ a_n\int_{T_y\wedge d\mu_{q-l}}\biggl( \int_{(T_x\wedge d\mu_{p-k})\bigl(\zeta(\mathbf x, \mathbf y)-
\zeta_x^n (\mathbf x)\zeta_y^n(\mathbf y)\bigr)
} \phi(\mathbf x_k, \mathbf y_l) d\mu_x\biggr) d\mu_y$$ for all $n$  is bounded. Thus the sequence of real numbers
$$\int_{T_y\wedge d\mu_{q-l}}\biggl( \int_{(T_x\wedge d\mu_{p-k})\bigl(\zeta(\mathbf x, \mathbf y)-
\zeta_x^n (\mathbf x)\zeta_y^n(\mathbf y)\bigr)
} \phi(\mathbf x_k, \mathbf y_l) d\mu_x\biggr) d\mu_y$$ converges to $0$ as $n\to \infty$.
Taking the limit as $n\to \infty$,  the formula  (2.67) becomes
\begin{equation}
\int_{\mathcal T_{k, l}(\xi)}\phi=\displaystyle{\lim_{n\to\infty}}\int_{(\mathbf x_k, \mathbf y)\in V_{k, l}} 
\mathcal L_k^n(\mathbf x_k)\mathcal L_l^n (\mathbf y_l)\phi(\mathbf x_k, \mathbf y_l) d\mu_k\wedge d\mu_l.\end{equation}
(the Lebesgue integral exists due to the part (1) of Lemma 2.13).   Then the inequality 
$$\int_{\mathcal T_{k, l}(\xi)}\phi\leq C||\phi||_{0, K}$$ for some constant $C$ holds. 
By Proposition 2.1.8, 1.3.11, [4], $\mathcal T_{k, l}(\xi)$ is a distribution  of order $0$, thus a measure. 
If $\phi$ is a characteristic function of a set with $0$ Lebesgue measure. By the Fubini's theorem, $$\int_{\mathcal T_{k, l}(\xi)}
\phi=0.$$
  Thus $\mathcal T_{k, l}(\xi)$ is absolutely continuous with respect to the Lebesgue measure. 
The Lebesgue integral expression
(2.70) also shows that the Radon-Nilkodym derivative has inequality
\begin{equation}
|{d\mathcal T_{k, l}(\xi)\over d\mu_{k, l}}|\leq    C'M^2
\end{equation}
for some constant $C'$. 
We complete the proof of the 
Lebesgue condition.

\bigskip

Next we prove  the Radon-Nikodym condition.  

Let  $$\mathcal L_{k, l}(\mathbf x_k, \mathbf y_l)={d\mathcal T_{k, l}(\xi)\over d\mu_{k, l}}$$
 be
the Lebesgue function, i.e. the projection 
of De Rham distributions of \begin{equation}
\biggl((\xi_xT_1)(\xi_y T_2)\biggr)\wedge \xi\end{equation}
to $V_{k, l}$.   
Let's consider the integral
\begin{equation}
A_{\bm\lambda_1, \bm\lambda_2}=\int_{V_{k, l}} \mathcal L_{k, l}(D_{\bm\lambda_1 }(\mathbf x_k), 
D_{\bm\lambda_2}(\mathbf y_l)) \phi(\mathbf x_k, \mathbf y_l) d\mu_{k, l}.
\end{equation}
for a test function $\phi$ on $V_{k, l}$, where
$$\bm\lambda_1=\lambda_1^1 \mathbf e_1+\cdots+\lambda_k^1 \mathbf e_k, \quad \lambda_i^1\in \mathbb R^+$$
$$\bm\lambda_2=\lambda_1^2 \mathbf b_1+\cdots+\lambda_l^2 \mathbf b_l, \quad \lambda_i^2\in \mathbb R^+.$$
and $D_{\bm\lambda_1}, D_{\bm\lambda_2}$ are testing  maps as in (2.3). 
By Theorem 6, Chapter II, \S 7, [3], there are sequences of test forms
\begin{equation}
\xi_x^n (\mathbf x), \xi_y^n(\mathbf y)
\end{equation}
on $\mathcal X, \mathcal Y$ respectively 
such that $$ \xi_x^n (\mathbf x) \xi_y^n(\mathbf y)\to \xi, \ as\ n\to \infty. $$
The projection of a De Rham distribution of 
\begin{equation}
\biggl((\xi_xT_1)(\xi_y T_2)\biggr)\wedge \biggl(\xi- \xi_x^n (\mathbf x) \xi_y^n(\mathbf x)\biggr)\end{equation}
to $V_{k, l}$ is 
an  $L^1$ function $\mathcal L^n (\mathbf x_k, \mathbf y_l)$.
Similarly as before, there is a sequence of number $a_n$ converges to $+\infty$ such that
$$a_n(\xi- \xi_x^n (\mathbf x) \xi_y^n(\mathbf x))$$ is bounded 
for all $\mathbf x, \mathbf y, n$.  By the Lebesgue condition of $\mathcal L^n (\mathbf x_k, \mathbf y_l)$, 
 $$a_n \mathcal L^n (\mathbf x_k, \mathbf y_l)$$ is bounded for all $n$. Hence 
 $$\mathcal L^n (\mathbf x_k, \mathbf y_l)$$ converges to $0$ uniformly. 
Hence \begin{equation}\begin{array}{c} B_{\bm\lambda_1, \bm\lambda_2}^n\\
\|\\
\int_{V_{k, l}}\mathcal L^n ( D_{\bm\lambda_1 }(\mathbf x_k), D_{\bm\lambda_2 }(\mathbf y_l))
\phi(\mathbf x_k, \mathbf y_l) d\mu_{k, l}.
 \end{array}\end{equation}
converges to zero as $n\to \infty$ uniformly for all $\bm\lambda_1, \bm\lambda_2$.

\bigskip

We let \begin{equation}\begin{array}{c} C^n_{\bm\lambda_1, \bm\lambda_2}\\
\|\\
 \int_{V_{k, l}} \mathcal L_x^n( D_{\bm\lambda_1}(\mathbf x_{k}))
\mathcal L_y^n(D_{\bm\lambda_2}(\mathbf y_l))\phi(\mathbf x_k, \mathbf y_{l})d\mu_{k, l}
\end{array}\end{equation}
where  $\mathcal L_x^n( \mathbf x_{k}), 
\mathcal L_y^n(\mathbf y_l)$ are the projections of De Rham distributions of 
$$\xi_x T_x\wedge \xi_x^n,  \xi_y T_y\wedge \xi_y^n$$
to $V_k, V_l$ respectively. 
By part (1) of Lemma 2.13, 
on the individual manifold of $\mathcal X, \mathcal Y$,  for each fixed $n$
there is an iterated limit
\begin{align*}
& \displaystyle{\lim_{|\bm\lambda_2|\Rsh 0}}
\int_{\mathbf y\in V_l}\mathcal L_y^n(D_{\bm\lambda_2}(\mathbf y_l))\\ &
\biggl\{\displaystyle{\lim_{|\bm\lambda_1|\Rsh 0}} \int_{\mathbf x\in V_k} 
\phi(\mathbf x_k, \mathbf y_l) \mathcal L_x^n( D_{\bm\lambda_1}(\mathbf x_{k})) d\mu_k\biggr\}
d\mu_{l}.
\end{align*}

By the uniform convergence of Lemma 2.13,  this limit can be written as a single limit  
\begin{equation}\begin{array}{c}
 R_n\\\|\\
\displaystyle{\lim_{|(\bm\lambda_1, \bm\lambda_2)|\Rsh 0}} C^n_{\bm\lambda_1, \bm\lambda_2}\\
\|\\
\displaystyle{\lim_{|(\bm\lambda_1, \bm\lambda_2)|\Rsh 0}}
\int_{V_{k, l}} \mathcal L_x^n(D_{\bm\lambda_1}(\mathbf x_{k}))
\mathcal L_y^n(D_{\bm\lambda_2}(\mathbf y_l))\phi(\mathbf x_k, \mathbf y_{l}) 
 d\mu_{k, l}.\end{array}\end{equation}
Let 
$$L=inf\displaystyle{\lim_{n\to \infty}}R_n,$$
which is finite (because $R_n$ is bounded). 

For a given $\epsilon>0$, there are infinitely many $n\geq N_1$ such that 
\begin{equation}
|R_{n}-L|\leq \epsilon.
\end{equation}

Hence there is a particular  $N_3$ such that 
both inequalities\begin{equation}\begin{array}{c}
|R_{N_3}-L|\leq \epsilon\\
|B_{\bm\lambda_1, \bm\lambda_2}^{N_3}|\leq \epsilon
\end{array}\end{equation}
hold for all $\bm\lambda_1, \bm\lambda_2$. 
Applying Lemma 2.13 with the fixed $N_3$, we can find points $(\bm\lambda_1, \bm\lambda_2)$ on the zigzag path with $|||(\bm\lambda_1, \bm\lambda_2)|||\prec\epsilon$,
\begin{equation}
 |C^{N_3}_{\bm\lambda_1, \bm\lambda_2}-R_{N_3}|\leq \epsilon.
\end{equation}

Then three inequalities (2.80), (2.81) connect four numbers, 
$$A_{\bm\lambda_1, \bm\lambda_2}\to  C^{N_3}_{\bm\lambda_1, \bm\lambda_2}\to  R_{N_3}\to  L$$
in the following conclusion:  for any positive $\epsilon$, there exist number $N_3$ followed by the 
points $(\bm\lambda_1, \bm\lambda_2)$ on zigzag path with 
$$|||(\bm\lambda_1, \bm\lambda_2)|||\prec\epsilon$$
such that 
\begin{equation}
|A_{\bm\lambda_1, \bm\lambda_2}-L|\leq 3\epsilon.\end{equation}

So 
$$
\displaystyle{\lim_{|(\bm\lambda_1, \bm\lambda_2)|\Rsh 0}} A_{\bm\lambda_1, \bm\lambda_2}=L.
$$
We complete the proof.

\end{proof}
\bigskip

\begin{proposition}
If $T$ is Lebesgue and $\omega$ is $C^\infty$,  then the intersection 
\begin{equation}
T\wedge \omega
\end{equation}
is Lebesgue.

\end{proposition}

\bigskip

\begin{proof}   This is the tautology. 
Let $\xi\in \mathscr D^s(U)$.
Notice $\omega\wedge \xi\in \mathscr D(U)$.
Then the projection $\mathcal J$ of De Rham distributions of 
$$(T\wedge \omega)\wedge \xi$$
is the same as that of 
$$ T\wedge (\omega \wedge \xi)$$
Since $\mathcal J$ satisfies both 
Lebesgue condition and Radon-Nikodym condition for $T$ is Lebesgue.   Hence $T\wedge \omega$ is Lebesgue. 
  We complete the proof. 

\end{proof}

\bigskip

\bigskip

\begin{ex} There exist  currents that are not Lebesgue.

In the Euclidean space $\mathbb R^m$ of coordinates $x_1, \cdots, x_p, \cdots, x_m$,  
we let $$T=\delta_O dx_{p+1}\wedge \cdots\wedge dx_m$$
with $\delta$-function $\delta_O$ of the origin $O$ of $\mathbb R^m$.  Let $V$ be the coordinates plane with coordinates 
$x_1, \cdots, x_p$, and $\pi: \mathbb R^m\to V$ be the projection.
 Let $\xi\in \mathscr D(\mathbb R^m)$ with  $\xi(O)\neq 0$.
So a projection of the De Rham distribution $\pi_\star (\xi \delta_O)$ is equal to
\begin{equation}
\delta_{\mathbf 0} \xi(O) 
\end{equation}
where $\delta_{\mathbf 0}$ is the $\delta$-function at the origin in the plane $V$. 
Hence $\pi_\star (\xi \delta_O)$ is the distribution $\delta_{\mathbf 0}\xi(O)$,   also a measure  on $V$ with the Borel $\sigma$-algebra of $V$.  
Now we consider the two measures for $V$ on the same $\sigma$-algebra.  When they are applied to the singleton set, the origin of $V$, the Lebesgue  measure  is zero, but the projection  measure  is $\xi(O)\neq 0$.  
Hence $$\pi_\star (\xi \delta_O)\not\ll\text{ Lebesgue  measure}.$$  So $T$ 
does not satisfy the Lebesgue condition. \end{ex}

\bigskip

\section { Regularization}\bigskip
 G. de Rham has  introduced the notion--current that  connects the singular chains and $C^\infty$ forms. 
 In his original work ([3]), there  were two important insights on currents that led to the De Rham's theorem: a chain is the limit of $C^\infty$ forms; the limit preserves 
 a homotopy  between a chain
and a $C^\infty$ form.  The construction to realize them  is the De Rham's 
regularization, which consists of  the construction of two operators:  $R_\epsilon$, $A_\epsilon$ ( 
see chapter III, [3])
 \footnote{De Rham's regularization is based on the Schwartz's regularization of distribution,   \S 4, chapter VI, [5].}.  
We are going to  study the associated geometric measures. 
 Let's recall it.

\bigskip

\subsection{De Rham's construction}

\begin{definition} 
Let $\mathcal X$ be a connected, oriented manifold. 
Let $\epsilon$ be a small positive number.  Linear operators $R_\epsilon$ and   $A_\epsilon$ on$\mathscr D'(\mathcal X)$
are called De Rham's regulator  and homotopy operators respectively  if they satisfy
\par
(1) a homotopy formula
\begin{equation}
R_\epsilon T-T=b A_\epsilon T+ A_\epsilon bT.
\end{equation}
\par \hspace{1cc} where $b$ is the boundary operator. 
\par
(2) $supp(R_\epsilon T), supp(A_\epsilon T) $ are contained in any given neighborhood of\par \hspace{1cc}
$supp(T)$ provided $\epsilon$ is sufficiently small.\par
 
(3) $R_\epsilon T$ is $C^\infty$;\par
(4)  If $T$ is $C^r$, $A_\epsilon T$ is $C^r$.\par
(5) If a smooth differential form $\phi$ varies in a bounded set and $\epsilon$ is bounded \par \hspace{1cc} above, then
$R_\epsilon \phi, A_\epsilon \phi$ are bounded.\par
(6) as $\epsilon\to 0$, 
$$\int_{R_\epsilon T}\phi \to \int_{T}\phi, \int_{A_\epsilon T}\phi\to 0$$
 \par \hspace{1cc} uniformly on each bounded set $\phi$, where the integral symbol denotes  \par \hspace{1cc} 
the functional of the current.  
\end{definition}
\bigskip

\bigskip

\begin{theorem} (G. de Rham)
The operators $R_\epsilon, A_\epsilon$ exist.
\end{theorem}

\bigskip

\begin{proof} In the following we review the constructions of  operators $R_\epsilon$ and $A_\epsilon$.  
The verification of conditions (1)-(6) in [3] will be omitted.      
There are three steps in the construction.\par
Step 1: Local construction. Construction in $\mathcal X=\mathbb R^m$. \par
Step 2:  Preparation. To prepare for the gluing process, convert the \par \hspace{2CC} \quad  operators to  a bounded domain $B$ in $\mathbb R^m$ with boundary. 

 \par
Step 3: Gluing. Assume $\mathcal X$ is covered by the bounded domain with  \par \hspace{2CC}\quad  boundary $B^i, \ countable\ i$.  Then  glue the operators in each $B^i$  to  \par \hspace{2CC}\quad  obtain 
the global  
\begin{equation}
R_\epsilon, A_\epsilon
\end{equation}

\bigskip

Step 1:  The most part of this step is originated from Schwartz's work in [5].   So we should've  skipped it. But
this is also the main technical base of this paper, and some details must be explored further, 
thus we are going to adapt it in a slightly different way.
  Let $\mathcal X=\mathbb R^m$ be the Euclidean space of dimension $m$ with the standard linear structure.  Let 
$x=(x_1, \cdots, x_m)$ be its
Euclidean coordinates, and vectors and points in $\mathbb R^m$  will be denoted by the $\bf bold$ letters.
 Let  $T$ be a homogeneous
current of degree $p$ on $\mathbb R^m$.  Let $ B$ be a bounded domain
in $\mathbb R^m$ with boundary.
 Let $\phi(\mathbf x)\in \mathscr D(B)$   satisfy
\begin{equation}
\int_{\mathbf x\in \mathbb R^m}\phi(\mathbf x)d\mu_x=1,
\end{equation} where
$d\mu_x$ is the Euclidean volume form 
$$dx_1\wedge\cdots \wedge dx_m.$$

Let $$f^{\epsilon}=\epsilon^{-m} f(\epsilon^{-1}\mathbf x),\quad \quad  \epsilon\in \mathbb R^+ $$

 Let 
\begin{equation}
\vartheta_\epsilon(\mathbf x)=f^{\epsilon}(\mathbf x)d\mu_x=\vartheta_1({\mathbf x\over \epsilon}).\end{equation} be the $m$-form on $\mathbb R^m$.
\bigskip

Next we define two operators on the differential forms of Euclidean space $\mathbb R^m$ based on  $C^\infty$ maps $s_{\mathbf y}(\mathbf x)$ below. 
 Let $$s_{\mathbf y}(\mathbf x)$$ be   $C^\infty$ maps parametrized by $\mathbf y\in \mathbb R^m$, 
$$\begin{array}{ccc}
\mathbb R^m &\rightarrow & \mathbb R^m\\
\mathbf x   &\rightarrow & s_{\mathbf y}(\mathbf x)\end{array}$$ such that
all partial derivatives of the components with respect to the variables of $\mathbf x$ are continuous functions
in $(\mathbf x, \mathbf y)$.   
Let $\phi$ be a test form on $\mathbb R^m$. 
For such maps $s_{\mathbf y}(\mathbf x)$, we denote two operations on the form $\phi$  
$$s_{\mathbf y}^\ast(\phi), \quad\quad \ and $$
$$\mathbf S_{\mathbf y}^\ast(\phi)=Proj_\ast ( s_{(t, \mathbf  y)}^\ast (\phi)), t\in [0, 1]$$
where $Proj: [0, 1]\times \mathcal X\to \mathcal X$ is the projection, and
$$\begin{array}{ccc}
 s_{(t, \mathbf y)}: [0, 1]\times \mathcal X &\rightarrow &\mathcal X\\
(t, x) &\rightarrow & s_{t\mathbf y}(\mathbf x).\end{array}$$

\bigskip

 Then we define  operators $R_\epsilon, A_\epsilon$ on currents $T$ by

\begin{equation}\begin{cases} \int_{R_\epsilon   T}\phi=
\int_{\mathbf x\in T} \biggl(\int_{\mathbf y\in \mathbb R^m} \vartheta_\epsilon(\mathbf y)
\wedge s_{\mathbf y}^\ast \phi(\mathbf x)\biggr),\\
\int_{A_\epsilon T}\phi =\int_{\mathbf x\in T} \biggl( \int_{\mathbf y\in \mathbb R^m}  
\vartheta_\epsilon(\mathbf y)\wedge \mathbf S_{\mathbf y}^\ast\phi(\mathbf x)\biggr)
\end{cases}\end{equation}
where $\phi$ is a test form. 
We should note that \par
(1) the continuity assumption about $s_{\mathbf y}(\mathbf x)$ guarantees the existence of \par\hspace{1CC}  (3.5),\par
(2)  also equations \begin{equation}\begin{cases}
deg(s_{\mathbf y}^\ast(\phi))=deg(\phi),\\
deg(\mathbf S_{\mathbf y}^\ast(\phi))=deg(\phi)-1\end{cases}
\end{equation}

\noindent imply that 
\begin{equation}\begin{cases}
dim (R_\epsilon(T))=dim(T), \\
 dim(A_\epsilon(T))=dim(T)-1.\end{cases}
\end{equation}
If furthermore the map
$$\begin{array}{ccc}
\mathbb R^m\times \mathbb R^m &\rightarrow &  \mathbb R^m \times \mathbb R^m\\
(\mathbf x, \mathbf y)   &\rightarrow & (\mathbf x, s_{\mathbf y}(\mathbf x))\end{array}$$
is a diffeomorphism,  there is a  change of variables
\begin{equation}\begin{cases}
s_{\mathbf y}(\mathbf x)\Rightarrow \mathbf x\\
\mathbf  y\Rightarrow s^{-1}(\mathbf x, \mathbf y)
\end{cases}\end{equation}
(replacement of $s_{\mathbf y}(\mathbf x)$ with $\mathbf x$; $\mathbf  y$ with $s^{-1}(\mathbf x, \mathbf y)$)
where $s^{-1}:\mathbb R^m\times \mathbb R^m\to \mathbb R^m$ is $C^\infty$ and 
satisfies $s_{s^{-1}(\mathbf x, \mathbf y)}(\mathbf x)=\mathbf y$. Then 
 the first integral of (3.5) shows that 
\begin{equation} {R_\epsilon   T}=\int_{\mathbf x\in T}  \vartheta_\epsilon(s^{-1}(\mathbf x, \mathbf y))
\end{equation}
 is a $C^\infty$ form (see p65, [3]).  The form $\vartheta_\epsilon(s^{-1}(\mathbf x, \mathbf y))$ as a form in variables  
 $\mathbf x, \mathbf y$  is 
 the kernel (p71, [3]) of $R_\epsilon$. 
We should make a note that the current's evaluation  (3.9)  is defined in the same way as the fibre integrals of $C^\infty$ forms  under the projection $\mathcal X\times \mathcal X\to \mathcal X$.

\bigskip

While  a general $s_{\mathbf y}(\mathbf x)$ will be used in later step, in the step 1 
we  use  $s_{\mathbf y}(\mathbf x)=\mathbf x+\mathbf y$ for the particular case of $\mathbb R^m$, 
where the $+$ is from the standard linear structure of $\mathbb R^m$. Then $R_\epsilon$ is  the convolution. 
Next we sketch the rest of two steps in the globalization.\bigskip

Step 2:  Choose the unit ball $B\subset \mathbb R^m$ diffeomorphic to $\mathbb R^m$. Let $h$ be the specific diffeomorphism  $$
\begin{array}{ccc}
\mathbb R^m &\rightarrow &  B,\end{array}$$
defined on p66, [3].
Denote the $s_{\mathbf y}(\mathbf x )$ in step 1 by $s_{\mathbf y}^+(\mathbf x )$. Then we define the new $C^\infty$ map
 
\begin{equation}s_{\mathbf y}(\mathbf x)=\left\{ \begin{array}{ccc} h s_{\mathbf y}^+ h^{-1} (\mathbf x) & \mbox{for } &   \mathbf x\in B\\
\mathbf x & \mbox{for }  & \mathbf x\notin B.
\end{array}\right.
\end{equation}
We would like to point out that  $s_{\mathbf y}(\mathbf x)$ satisfies  assumption.  
Then we can define the operators $R_\epsilon^B, A_\epsilon^B$ depending on $B$ in the same way (with a test form $\phi$):

\begin{equation}\begin{cases} \int_{R_\epsilon^B   T}\phi=
\int_{\mathbf x\in T} \biggl(\int_{\mathbf y\in \mathbb R^m} \vartheta_\epsilon(\mathbf y)
\wedge s_{\mathbf y}^\ast \phi(\mathbf x)\biggr),\\
\int_{A_\epsilon^B T}\phi =\int_{\mathbf x\in T} \biggl( \int_{\mathbf y\in \mathbb R^m}  
 \vartheta_\epsilon(\mathbf y)\wedge \mathbf S_{\mathbf y}^\ast\phi(\mathbf x)\biggr).
\end{cases}\end{equation}

Then the operators $R_\epsilon^B, A_\epsilon^B$  will satisfy\par

(a) properties  (1), (4), (5) and (6) in definition 3.1.\par

(b) $ R^B_\epsilon(T)$ is $C^\infty$ in $B$, $R^B_\epsilon(T)=T$ in the complement of
$\bar B$;\par
(c)  if $T$ is $C^r$ in a neighborhood of a boundary point of $B$, $A^B_\epsilon(T)$ will have the same regularity 
in the neighborhood. 
(There is a slight  difference in notations from [3]).

\bigskip

Step 3: Cover the  $\mathcal X$ with   countable open sets $B_i$ (locally finite). Now we regard each  $B^i$ as a subset of $B$ in step 2.
Let a neighborhood $U_i$ of $B_i$.  Let $h_i$ be the  diffeomorphic-to-image map
$$\begin{array}{ccc}
U_i &\rightarrow  &\mathbb R^m\\
\cup   && \cup\\
B_i &\rightarrow & B.
\end{array}$$
Let $g_i\geq 0$ be a function on $\mathcal X$, which is 1 on  $B_i$ and supported in $U_i$.  Let
$T'=g_iT$ and $T''=T-T'$.
Then we let 
$$\begin{array}{c}R^i_\epsilon T=(h_i^{-1})_\ast \circ R_\epsilon^B  \circ (h_i)_\ast T'+T''\\

A^i_\epsilon T=(h_i^{-1})_\ast  \circ A_\epsilon^B  \circ (h_i )_\ast T'.\end{array}$$
( Note: $h_i^{-1}$ is well-defined because $h_i$ is a diffeomorphic-to-image map). 
Finally we glue them together by taking the composition,

$$ \begin{array}{c}
R^{(\mathcal N)}_\epsilon=R^1_{\epsilon}\circ \cdots \circ R^{\mathcal N}_{\epsilon}, \\
A^{(\mathcal N)}_\epsilon=R^1_{\epsilon}\circ \cdots \circ R^{\mathcal N}_{\epsilon} \circ A^\mathcal N_\epsilon.
\end{array}
$$
Then we take the  limit as $\mathcal N\to \infty$ with respect to the compact support to obtain the well-defined, global
operator $R_\epsilon$ and $A_\epsilon$.\footnote{ In [3], for each open set $U_i$ there is a different positive $\epsilon_i$. 
We used the same number $\epsilon$ for all $U_i$. This difference should be noticed.}   

\smallskip

 \end{proof}

\bigskip

\begin{definition} ( De Rahm data for De Rham's regularization)\quad\par
(a) We call $R_\epsilon$ from Theorem 3.2 the De Rham's regulator,  \par\hspace{1cc} 
$ A_\epsilon$ from Theorem 3.2 the De Rham's homotopy operator,   
 and the \par\hspace{1cc} 
associated regularization the De Rham's   regularization. All operators \par\hspace{1cc} 
$R_\epsilon, A_\epsilon$ in this paper are chosen to be De Rham's. (The general  \par\hspace{1cc} operators 
from definition 3.1 are not necessarily De Rham's).  \par
(b)    
We define De Rham data to be all items in the construction of  \par\hspace{1cc} 
De Rham's regularization operators  $R_\epsilon, A_\epsilon$.  More specifically it  \par\hspace{1cc} includes   \par\hspace{1cc}
(1) the covering $B_i\subset U_i$ with the order of countable $i$,  \par\hspace{1cc}
(2) the diffeomorphism $h_i: U_i\to \mathbb R^m$, and  functions $g_i$ with value  \par\hspace{1cc}\quad \hspace {1.5pt}
 1  on $B_i$,  \par\hspace{1cc}
(3) for each $B_i$, another diffeomorphism $h^i: B^i\simeq \mathbb R^m$ with  Euclidean 
\par\hspace{1cc}\quad \hspace {1.5pt}  coordinates.  \par\hspace{1cc}
(4)  functions $f_i$ in each $B^i\simeq \mathbb R^m$ as in the first step, called  \par
  \par\hspace{1cc}\quad \hspace {1.5pt}  convolution functions, function $g_i$ called  the gluing functions. \par
(c) the covering $B_i\subset U_i$ equipped with all the items (1)-(4) in De Rham \par\hspace{1cc} \hspace {1pt}  
data is called the De Rham covering. Each pair $B_i\subset U_i$ with (1)-(4)   \par\hspace{1cc} \hspace {1pt}is called a   
De Rham chart.

\end{definition}
 
\bigskip

{\bf Remark} 
 The De Rham data gives  a covering that regularizes  the piece of $T$ supported in $B_i$ in each chart $U_i$ independently and leave 
$T$ outside of $B_i$ untouched.  
There is ``glue" (such as $g_i$) at each chart to piece them together.  But there is no 
relation among the charts. This extrinsic nature  is contrary to the other cases where a condition of the invariance must be imposed
on each chart to keep it intrinsic.   
 \bigskip

G. de Rham further showed in chapter III, \S 17, [3], \bigskip

\begin{corollary}
The De Rham's operator $R_\epsilon$ constructed in  Theorem 3.2 is a regularizing operator, i.e.
there is a $C^\infty$ form $\varrho_\epsilon(\mathbf x, \mathbf y) $ on $\mathcal X\times \mathcal X$, 
called the $C^\infty$ kernel of $R_\epsilon$,  such that as currents, 
$$R_\epsilon T =\int_{\mathbf y\in T} \varrho_\epsilon(\mathbf x, \mathbf y).$$

\end{corollary}

\bigskip

{\bf Remark} Kernel of an operator is a crucial technical notion defined by De Rham in [3] 
who focused on Hodge's harmonic integral theory.
The notion is also crucial in our new direction. So we list its definition in the Appendix.

\bigskip

\subsection{Kernel of De Rham's regulator}
\bigskip

\begin{definition}
Let $\omega$ be a $C^\infty$ form of degree $p$ on a manifold $\mathcal X$.
We say $\omega$ is a local constant slicing, if at each point, there is  a chart  $U$ containing the point such that 
\begin{equation}
\omega|_U=\sum_{finite \ k} \pi_k^\ast \biggl(\theta_k(\mathbf v_k)|_{\pi_k(U)}\biggr)
\end{equation}
where $\pi_k: U\to  V_k\simeq \mathbb R^p$ is the projections to the coordinates planes $V_k$ of dimension $p$, and 
$\theta_k(\mathbf v_k)$ is a $C^\infty$ form on $V_k$ (with points $ \mathbf v_k\in V_k$).
\end{definition}
\bigskip

{\bf Remark}  A form of a local constant slicing is a particular type of forms invariant under the $C^\infty$ diffeomorphisms.   
 For instance  we  notice that the differential operation commutes with the projection map $\pi_I$: 
$$d \omega|_U=\sum_{finite \ k} \pi_k^\ast \biggl( d\theta_k(\mathbf v_k)|_{\pi_k(U)}\biggr)$$
and $\theta_k(\mathbf v_k)$ is a form of maximal degree on $V_k$. Therefore 
a form of a local constant slicing must be closed. Hence it represents a cohomology class. 
\bigskip

\begin{lemma}
Let $\mathcal X_0$ be the union of countably many proper submanifolds of dimension strictly less than $dim(\mathcal X)$.  
Let  $\omega$ be a  $C^\infty$ form on $\mathcal X$  such that for each point $q\in \mathcal X$  there is a chart $U$ containing $q$ and 
all forms in (3.12) are well-defined $C^\infty$ forms on $U$, but 
the equality (3.12) only holds 
on the submanifold $U-(\mathcal X_0\cap U)$.   Then $\omega$ is still a local constant slicing on $\mathcal X$.
 \end{lemma}
\bigskip

\begin{proof} Let $U$ be the neighborhood in a chart. 
Then the restriction to $U-U\cap \mathcal X_0$ is the same  chart on the subset 
$U-\mathcal X_0\cap U$. 
By the assumption there are $C^\infty$ forms $\theta_k$ of maximal degree on coordinates planes $V_k$  such that
\begin{equation}
\omega|_{U-U\cap \mathcal X_0}=\sum_{finite \ k} \pi_k^\ast \biggl (\theta_k(\mathbf v_k)|_{\pi_k(U-U\cap \mathcal X_0)}\biggr)
\end{equation}
Notice both sides have extension to $U$ by the continuity.  Taking the  closure (of topology of $\mathcal X$) both sides, we complete the proof. 
\end{proof}

\bigskip

 We'll show the kernel of $R_\epsilon$ is not only $C^\infty$, but also a  local constant slicing, therefore closed. 
\bigskip

\begin{proposition}
Let $\varrho_\epsilon$ be the $C^\infty$ kernel of De Rham's regulator $R_\epsilon$. Then
$\varrho_\epsilon$ is  a  local constant slicing. Furthermore there are two charts $x, y$ in the De Rham data for a neighborhood $U_i$ at each point such that 
\begin{equation}\varrho_\epsilon|_{U}(\mathbf x, \mathbf y)=\varrho_1|_{U}({\mathbf x\over \epsilon}, {\mathbf y\over \epsilon})\end{equation}
inside of $B_i$ (where $x, y$ charts could the same chart).

\end{proposition}
\bigskip

{\bf Remark}. The $C^\infty$ kernel $\varrho_\epsilon$ is a closed form.  By the homotopy formula (3.1) it  represents the class of
the diagonal in the cohomology group with compact support.  
\bigskip

\begin{proof}
 For this particular local constant slicing $\varrho_\epsilon$, we'll give a full detailed description of the local projection for the constant slicing. \bigskip

Denote the boundary of each local ball $B_i$ in the De Rham data by $\partial _i$.  Let $\partial =\sum_i \partial_i$.
  By Lemma 3.6, it suffices to consider the submanifold $\mathcal X-\partial$.
So let $q\in \mathcal X-\partial$. 
Let $U_q\subset \mathcal X-\partial$ be a small neighborhood of $q$.   Consider the 
kernel $\varrho_\epsilon^{q}(\mathbf x,
\mathbf y)$
of the De Rham's regulator \begin{equation}
R_\epsilon=R_\epsilon^{1}\circ \cdots \circ R_\epsilon^{\mathcal N}\end{equation}
 restricted to $U_q\times U_q$, where $\mathcal N$ is finite because the covering is locally finite. 
Because we exclude $\partial$, there are two cases for the points $q$. If $q\not \in B_i$ for some $i$, 
$R_\epsilon^{i}|_{U_q}$ by the definition is the identity. If $q\in B_i$ for some $i$, 
$R_\epsilon^{i}|_{U_q}$,  which will be called regulator, has the $C^\infty$ 
kernel $\varrho^i_\epsilon(\mathbf x, \mathbf y)$ where 
$\mathbf y$ is in the second copy of $U_q$.  Suppose there are $n$ regulators in (3.15), and they
are in $B_1, B_2, \cdots, B_n$.
Let's denote the coordinates for each $U_i\supset B_i$ by the same letter $\mathbf x_i$  for which we should restrict ourselves to the domain $B_i$ only i.e.
apply the (3.10) for $\mathbf x_i\in B_i$ only. 
The kernel of each $R_\epsilon^{i }$ can be denoted by 
$$\vartheta_1^i({\mathbf x_i\over \epsilon}\overset{i}{-}{\mathbf y_i\over \epsilon}), $$ 
which says that for a current $T$,
$${R_\epsilon^{i}   T}=\int_{\mathbf y_i\in T}  
\vartheta_1^i({\mathbf x_i\over \epsilon}\overset{i}{-}{\mathbf y_i\over \epsilon})
$$
where the subtraction $\overset{i}{-}$ (  also  $\overset{i}{+}$), scalar multiplication  ${\bullet \over \epsilon}$ 
are from the linear structure of $U_i$ in De Rham data (they are come from the De Rham data). 

Hence the kernel $\varrho_\epsilon$ of $R_\epsilon=R_\epsilon^{1}\circ\cdots\circ R_\epsilon^{\mathcal N}$ inside $B_1\cap \cdots\cap B_n$
is the fibre integral
\begin{multline}
\varrho_\epsilon
=
 \int_{(\mathbf x_2, \cdots, \mathbf x_{n})\in \mathbb (R^m)^{\oplus n-1} } 
\vartheta_1^1({\mathbf x_1\over \epsilon}\overset{1}{-}{\mathbf x_{2}\over\epsilon})\wedge
\vartheta_1^{2}({\mathbf x_{2}\over\epsilon}
\overset{2}{-}{\mathbf x_{3}\over\epsilon})\wedge 
\cdots \\ 
\wedge \vartheta_1^{n-1} ({\mathbf x_{n-1}\over \epsilon}
\overset{n-1}{-}{\mathbf x_{n}\over \epsilon})
\wedge
\vartheta_1^n({\mathbf x_n\over \epsilon}\overset{n}{-}{ \mathbf y_{n}\over\epsilon}),
\end{multline}
whose degree is $m$. 
So $\varrho_\epsilon$ is 
the fibre integral of the local $C^\infty$ form,   
$$\begin{array}{c}\vartheta_1^1({\mathbf x_1\over \epsilon}\overset{1}{-}{\mathbf x_{2}\over\epsilon})\wedge
\vartheta_1^{2}({\mathbf x_{2}\over\epsilon}
\overset{2}{-}{\mathbf x_{3}\over\epsilon})\wedge 
\cdots 
\wedge \vartheta_1^{n-1} ({\mathbf x_{n-1}\over \epsilon}
\overset{n-1}{-}{\mathbf x_{n}\over \epsilon})
\wedge
\vartheta_1^n({\mathbf x_n\over \epsilon}\overset{n}{-}{ \mathbf y_{n}\over\epsilon})\\
\|\\
\vartheta_\epsilon^1({\mathbf x_1}\overset{1}{-}{\mathbf x_{2}})\wedge
\vartheta_\epsilon^{2}({\mathbf x_{2}}
\overset{2}{-}{\mathbf x_{3}})\wedge 
\cdots 
\wedge \vartheta_\epsilon^{n-1} ({\mathbf x_{n-1}}
\overset{n-1}{-}{\mathbf x_{n}})
\wedge
\vartheta_\epsilon^n({\mathbf x_n}\overset{n}{-}{ \mathbf y_{n}})
\end{array}$$
denoted by
\begin{equation} \varsigma_\epsilon^{(q)},\end{equation}
 (of degree $mn$), 
in the projection  of the Cartesian product
\begin{equation} \mathcal P_1: (\mathbb R^m)^{\oplus ( n+1)}\to \mathbb R^m_{\mathbf x_1}\oplus \mathbb R^m_{\mathbf y_n}\end{equation}
where  $(\mathbb R^m)^{\oplus ( n+1)}$ have  global coordinates 
$$\mathbf x_1, \cdots, \mathbf x_n, \mathbf y_n,$$ 
and $\mathbb R^m_{\mathbf x_1}, \mathbb R^m_{\mathbf y_n}$ are the first and last copies. 
Above argument is a technical description  of the kernel $\varrho_\epsilon$.   \bigskip

 To  associated a local constant slicing form,  we construct a commutative diagram by first  defining the diffeomorphism
$$\begin{array}{ccc}
\kappa_1: (\mathbb R^m)^{\oplus n+1} &\rightarrow &   (\mathbb R^m)^{\oplus n}\oplus \mathbb R^m\\
(\mathbf x_1, \cdots, \mathbf x_n, \mathbf y_n) &\rightarrow & ({\mathbf x_1}\overset{1}{-}{\mathbf x_{2}}, 
\cdots, {\mathbf x_n}\overset{n}{-}{ \mathbf y_{n}}, \mathbf y_n),
\end{array}
$$
where $\mathbf y_n$ are the coordinates for the last copy $\mathbb R^m$, and each copy $\mathbb R^m$ has its own 
linear structure. 
Then the formula (3.16) is just
$$\varrho_\epsilon
=(\mathcal P_1)_\ast (\varsigma_\epsilon^{(q)}).$$

We denote the components
in the target space $(\mathbb R^m)^{\oplus n}\oplus \mathbb R^m$ by
$$\mathbf x_1', \cdots, \mathbf x_n', \mathbf y_n.$$
Notice the map has rank $m(n-1)$, and
$\varsigma_\epsilon^{(q)}$ is the pullback form by $\kappa_1$: 

\begin{equation}
\varsigma_\epsilon^{(q)}=\vartheta_\epsilon^1(\mathbf x'_1)\wedge
\vartheta_\epsilon^{2}({\mathbf x'_{2}})\wedge 
\cdots 
\wedge
\vartheta_\epsilon^n({\mathbf x'_n}).
\end{equation}

So there is a commutative diagram

\begin{equation}\begin{array}{ccc}
(\mathbb R^m)^{\oplus n+1} &\stackrel{\kappa_1} = & (\mathbb R^m)^{\oplus n} \oplus \mathbb R^m\\
\scriptstyle{\mathcal P_1} \downarrow & &\scriptstyle{( \mathcal P_2, id)} \downarrow\\
\mathbb R^m_{\mathbf x_1}\oplus \mathbb R^m_{\mathbf y_n} &\stackrel{(\kappa_2, id)}=  
&\mathbb R^m \oplus \mathbb R^m_{\mathbf y_n} 
\end{array}\end{equation}
where $$\kappa_2: (\mathbf x_1, \mathbf y_n)\to {\mathbf x_1}\overset{n}{-}{ \mathbf y_{n}}$$ and

$$\mathcal P_2: (\mathbf x'_1, \cdots, \mathbf x'_n)\to \mathbf x'_1\overset{1}{+}\mathbf x'_2\overset{2}{+}
\cdots \overset{n-1}{+}\mathbf x'_n$$ is the map onto the first copy $\mathbb R^m$.
Then the commutativity of (3.20) yields \begin{equation}
\varrho_\epsilon=(\mathcal P_1)_\ast (\varsigma_\epsilon^{(q)})=( \mathcal P_2, id)_\ast (\varsigma_\epsilon^{(q)})\end{equation}
( the second equality in (3.21) is the isomorphism  
$$(\kappa_2, id): \mathbb R^m_{\mathbf x_1}\oplus \mathbb R^m_{\mathbf y_n}\to 
\mathbb R^m \oplus \mathbb R^m_{\mathbf y_n}. )$$ 
Since left hand side is a local constant slicing, i.e. the pullback of the form from $\mathbb R^m$, so is the right hand side. 
This completes the proof.

\end{proof}
\bigskip

\bigskip

\section{Intersection}

\bigskip

\subsection{Convergence of regularization}

\bigskip

\begin{theorem} \quad\par
Let $\mathcal X$ be a manifold equipped with De Rham data.  Let  $T_1, T_2$ 
be two homogeneous Lebesgue currents of dimensions $p, q$ respectively.  \par
(1)  Let $\phi$ be a test form of degree $p+q-m$. Then\par

\begin{equation}
\displaystyle{\lim_{\epsilon\to 0}}\int_{T_1} R_\epsilon T_2\wedge \phi
\end{equation}
\par\hspace{1 CC} exists. \par

(2) If  a sequence of test forms $\phi_n, n\in \mathbb N$ with the same compact support \par\hspace{1 CC} 
converges to $0$, so does the sequence  of numbers,
\begin{equation}
\displaystyle{\lim_{\epsilon\to 0}}\int_{T_1} R_\epsilon T_2\wedge \phi_n
\end{equation}
\par\hspace{1 CC} as $n\to \infty$.

\end{theorem}
\bigskip

\begin{proof} \par

(1)  Let $T_1, T_2$ are homogeneous currents of dimensions $p,  q$ respectively.
Then
\begin{equation}
\int_{T_1} R_\epsilon T_2\wedge \phi=(-1)^m \int_{T_1T_2\wedge \phi}\varrho_\epsilon(\mathbf x, \mathbf y).
\end{equation}
By Proposition 3.7, the kernel $\varrho_\epsilon(\mathbf x, \mathbf y)$ of $R_\epsilon$  is a  local constant slicing.  Thus there exists countable,  locally finite open covering $U_i$ of $\mathcal X$  such that 
\begin{equation}
\varrho_\epsilon(\mathbf x, \mathbf y)|_{U_i\times U_i}=
\varrho_1({\mathbf x\over \epsilon}, {\mathbf y\over \epsilon})|_{U_i\times U_i}
=\sum_{finite \ k} \pi_k^\ast (\theta_k({\mathbf v_k\over \epsilon}))
,\end{equation}
where $\pi_k: U_i\times U_i\to V_k$ is the projection to some coordinates planes (see Definition 3.5 for the complete notations.)

Let  $p_\ell(\mathbf x)$ be a partition of unity for the open covering $U_i$.
Then $T_1=\sum_{\ell}p_\ell T_1$ and each $p_\ell T_1$ is a current  supported on the open set $U_i $. 
  Since $supp(\phi)$ is compact, there will be finitely many such $\ell$ contributing to the integral (4.1). Therefore it suffices to show 
 the convergence of one integral 
\begin{equation}
\int_{p_\ell T_1} R_\epsilon (T_2)\wedge \phi.
\end{equation}
 for a fixed $\ell$. 
 We denote $p_\ell T_1$ by $S$.  The formula (4.5) assures us that we may assume $T_2$ also has a compact 
support in $U$.
 By formula (4.4), it suffices to consider the number
\begin{equation}
\int_{S T_2\wedge\phi} \pi_k^\ast (\theta_k({\mathbf v_k\over \epsilon})).
\end{equation}
By Proposition 2.6, 2.11, 2.14, $S T_2\wedge\phi$ is Lebesgue in $U_i\times U_i$. Let
$\mathcal L_k (\mathbf v_k)$ be a Lebesgue function on $V_k$. Then (4.6) is equal to
\begin{align}
 &\int_{\mathbf v_k\in V_k} \mathcal L_k (\mathbf v_k)\theta_k  ({\mathbf v_k\over \epsilon})\\
&= \int_{\mathbf v_k\in V_k} \mathcal L_k (\epsilon \mathbf v_k) \theta_k ({\mathbf v_k})
\end{align}
By the Radon-Nikodym condition for $S T_2\wedge\phi$, the limit
\begin{equation}
 \displaystyle{\lim_{\epsilon\to 0}}\int_{\mathbf v_k\in V_k} \mathcal L_k (\epsilon \mathbf v_k) \theta(\mathbf v_k)
\end{equation} exists.
Taking the sum of it over finite $k$ we complete the proof of part (1).

(2) Now we replace the test form $\phi$ by a sequence  of test forms $\phi_n\to 0$ as $n\to \infty$, and all have
 the same compact support.   
So there is a sequence of real numbers $a_n$ with $\displaystyle{\lim_{n\to\infty}} a_n=+\infty$ and
$\{a_n\phi_n\}$ for $n\in \mathbb N$ is a set of forms locally bounded to order $0$ (p38, [3]). 
By part (1), 
\begin{equation}
\displaystyle{\lim_{\epsilon\to 0}}
\int_{T_1}R_\epsilon(T_2)\wedge a_n\phi_n= a_n \sum_k 
\displaystyle{\lim_{\epsilon\to 0}}\int_{\mathbf v_k\in V_k} \mathcal L_k^n (\epsilon \mathbf v_k) \theta_k({\mathbf v_k\over\epsilon})
\end{equation}
where $\mathcal L_k^n ( \mathbf v_k) $ is the same Lebesgue function in (4.7) but replacing $\phi$ by $\phi_n$.
Then by the definition of the Lebesgue condition, $|\mathcal L_k^n (\epsilon \mathbf v_k) |$ is bounded by a number 
independent of $n$. Thus
$$ a_n \displaystyle{\lim_{\epsilon\to 0}}
\int_{T_1}R_\epsilon(T_2)\wedge \phi_n$$ is bounded by a number independent of $n$.
Hence $$
\displaystyle{\lim_{n\to \infty}}
\displaystyle{\lim_{\epsilon\to 0}}
\int_{T_1}R_\epsilon(T_2)\wedge \phi_n=0.
$$

We complete the proof. 
\bigskip

\end{proof}
\bigskip

\bigskip

\subsection{The intersection}

\bigskip

\begin{definition}  Let $T_1, T_2$ be  homogeneous Lebesgue currents on a manifold $\mathcal X$ equipped with  De Rham data.  
  By Theorem 4.1, the functional on $\mathscr D(\mathcal X)$, 
$$\phi\to \displaystyle{\lim_{\epsilon\to 0}}
\int_{T_1}  R_{\epsilon}T_2 \wedge \phi$$ is  linear, continuous. Therefore we  define
the intersection current
\begin{equation}
[T_1\wedge T_2]
\end{equation}
by the formula
\begin{equation}
\int_{[T_1\wedge T_2]} \phi=\displaystyle{\lim_{\epsilon\to 0}}
\int_{T_1}  R_{\epsilon}T_2 \wedge \phi
\end{equation} 
for a test form $\phi$. 
 Hence there is a well-defined bilinear map
\begin{equation}\begin{array}{ccc}
\mathcal C(\mathcal X)\times \mathcal C(\mathcal X) &\rightarrow &  \mathcal D'(\mathcal X) \\
(T_1, T_2) &\rightarrow & [T_1\wedge T_2],
\end{array}\end{equation}
dependent of De Rham data, where  $\mathcal D'(\mathcal X)$ denotes the space of currents.

\end{definition}

\bigskip

{\bf Remark}
The intersection $[\cdot\wedge\cdot]$ and De Rham's regularization $R_\epsilon, A_\epsilon$ all depend
on the De Rham data in Definition 3.3. We'll omit the notation for this dependence by  fixing a data in general 
arguments, but will make a note in a particular case where the multiple De Rham data is necessary.  
\bigskip

\begin{proposition} If $T_1, T_2$ are Lebesgue, so is 
\begin{equation}
[T_1\wedge T_2].
\end{equation}

\end{proposition}

\bigskip

{\bf Remark} The proposition extends Proposition 2.15.
\bigskip

\begin{proof} 
We continue the setting in Theorem 4.1.  In particular let $U$ be one of the open  sets $U_i$ in the covering charts   that admits 
the  local constant slicing  for the kernel of De Rham's regulator  $R_\epsilon$. 
Let the currents $T_1, T_2$ have dimensions $p$, $q+r$ respectively.  Let 
 $\xi\in \mathscr D( U)$ have degree $r$.   
Let $K$ be a compact set in $U$. By Theorem 4.1, $[T_1\wedge T_2]$ is a current of dimension $p+q-m$.
By Definition 2.4 of Lebesgue current, we need to have a projection to 
a plane. So we let $\mathbb R^{p+q-m}$ be a coordinates plane of dimension $p+q-m$. 
Let $\phi\in \mathscr D (U)$ be a test form expressed 
as $$\phi=f(\mathbf v) d\mu_{s}$$
where $f(\mathbf v)\in \mathscr D(\mathbb R^{p+q-m})$ supported in $K$, and  $d\mu_s$ is the Euclidean volume form of an another arbitrary coordinates plane $V_s$ of dimension $p+q-m$. 
We denote $(\pi_I)^\ast (f(\mathbf v) )$ also by $f(\mathbf v)\in C^\infty (U)$.
Then we repeat the argument in Theorem 4.1. 
In particular  we may replace $T_1$ by current $S$ supported in $U$. 
Recall
$V_k$ are coordinates planes of $U$ for the  local constant slicing  of kernel $\varrho_\epsilon$
such that
\begin{equation}
\varrho_\epsilon(\mathbf x, \mathbf y)|_{U\times U}=
\varrho_1({\mathbf x\over \epsilon}, {\mathbf y\over \epsilon})|_{U\times U}
=\sum_{finite \ k} \pi_k^\ast (\theta_k({\mathbf v_k\over \epsilon}))
.\end{equation}(see (4.4)), where $\theta_k$ may be assumed to have a compact support. 
Then we applying the formula (2.18),  the projection of a De Rham distribution of 
$[(S\wedge\xi)\wedge T_2]$ to $\mathbb R^{p+q-m}$  is the functional
\begin{equation}\begin{array}{ccc}
f &\rightarrow &  \int_{[S\wedge T_2]\wedge \xi}f(\mathbf v) d\mu_s,\end{array}\end{equation}
 which is \begin{equation}\begin{array}{ccc}
f &\rightarrow & \displaystyle{\lim_{\epsilon\to 0}}\sum_k \int_{[ (S\wedge \xi) T_2]\wedge 
\theta_k } f(\mathbf v) d\mu_s.\end{array}\end{equation}
Next we choose such De Rham distributions $\mathcal T_J$ of $[(S\wedge \xi) T_2]$  that
the index $J$ satisfies 
$$d\mu^J\wedge \theta_\epsilon\wedge d\mu_s$$ has the volume form for the space
$$U\times U$$
where $d\mu^J$ are the Euclidean volume forms of $m$-dimensional coordinates planes
of De Rham distribution $\mathcal T_J$ ( there are finitely many such indexes $J$). 
Applying the Lebesgue condition for  $ (S\wedge \xi) T_2$, we can Let $$\mathcal L^k(\mathbf v_k, \mathbf v)$$ 
be  the Lebesgue function of 
$\sum_J \mathcal T_J$
on $V_k\times \mathbb R^{p+q-m}$. 
 Hence the functional (4.16) is 
\begin{equation}\begin{array}{ccc}
f &\rightarrow & \sum_k \displaystyle{\lim_{\epsilon\to 0}}\int_{(\mathbf v_k, \mathbf v)\in  V_k\oplus (\{0\}\oplus \mathbb R^{p+q-m} )} \mathcal L^k(\epsilon \mathbf v_k, \mathbf v)\theta_k (\mathbf v_k)\wedge  f(\mathbf v)  d\mu_{p+q-m}.\end{array}\end{equation}
Thus it suffice to consider  functional $\mathcal I_k$ for each $k$, 
\begin{equation}\begin{array}{ccc}
f &\rightarrow & \displaystyle{\lim_{\epsilon\to 0}} \int_{(\mathbf v_k, \mathbf v)\in  V_k\oplus (\{0\}\oplus \mathbb R^{p+q-m} )} \mathcal L^k(\epsilon \mathbf v_k, \mathbf v)\theta_k (\mathbf v_k)\wedge  f(\mathbf v)  d\mu_{p+q-m}.\end{array}\end{equation}
 Since $\mathcal L(\mathbf v_k, \mathbf v), \theta_k(\mathbf v_k)$ are bounded, 
\begin{equation}
\biggl|\mathcal I_k|_f\biggr|\leq C ||f||_{0, K}
\end{equation}
where $C$ is a constant, and  $||\bullet ||_{0, K}$ is the semi-norm.
The inequality  holds for all compact sets $K$ of $U$. Hence 
$\mathcal I_k$ is a distribution of order $0$. By Proposition 2.1.8, 1.3.11, [4], $\mathcal I_k$ is a measure. 
Let $E\subset \mathbb R^{p+q-m}$ be a set of Lebesgue measure zero and 
let $\chi_n$ be a sequence of smooth functions on $\mathbb R^{p+q-m}$ such that 
$0\leq \chi_n\leq 1$ and converges uniformly to the characteristic function $\chi(E)$ of $E$.
Let $B$ be a ball in $V_k$  such that $B\times \mathbb R^{p+q-m}$ contains 
$supp(\mathcal L_k(\mathbf v_k, \mathbf v))$.  Also denote the volume for of $V_k$ by $d\mu_k$  
Then we can estimate

\begin{align*}  &\biggl |\displaystyle{\lim_{\epsilon\to 0}}\int_{(\mathbf v_k, \mathbf v)\in  V_k\oplus (\{0\}\oplus \mathbb R^{p+q-m} )}
 \mathcal L^k(\epsilon \mathbf v_k, \mathbf v)\theta_k (\mathbf v_k)\wedge  \chi_n(\mathbf v)   d\mu_{p+q-m}\biggr|\\ \
\\&
\leq C \biggl |\int_{\mathbf v_k\in B, \mathbf v\in \mathbb R^{p+q-m}} \chi_n (\mathbf v)d\mu_k\wedge   d\mu_{p+q-m}\biggr|\\&
=C d\mu_k  |_B \biggl |\int_{\mathbf  x\in \mathbb R^{p+q-m}} \chi_n(\mathbf v)d\mu_{p+q-m}\biggr|\to 0, as\  n\ \to\infty \end{align*}
where $C$ is a constant. 
Therefore $\mathcal I_k$ for each $k$ is a measure absolutely continuous with respect to the Lebesgue measure. 
 Let's denote the Radon-Nikodym derivative
$${d \mathcal I_k\over d\mu_{p+q-m}}=\mathcal L_k^x.$$
By the above argument, 
\begin{equation}
\mathcal L_k^x(\mathbf v)=\displaystyle{\lim_{\epsilon\to 0}}\int_{\mathbf v_k\in  V_k }
 \mathcal L^k(\epsilon \mathbf v_k, \mathbf v)\theta_k(\mathbf v_k)
\end{equation}
Since  $\mathcal L^k(\mathbf v_k, \mathbf v), \theta_k (\mathbf v_k)$ are bounded, so is $\mathcal L_k^x$, 
So the current $[T_1\wedge T_2]$ satisfies the Lebesgue condition of Definition 2.4. 
\bigskip

Now we continue to Radon-Nikodym condition. 
Let $$\mathbf e_1, \cdots, \mathbf e_{m}$$ be a linear basis for the Euclidean space $V_k$, and
$$\mathfrak a_1, \cdots, \mathfrak a_{p+q-m}$$
be a basis for $\mathbb R^{p+q-m}$.
Let \begin{equation}\begin{array}{c}
\bm\lambda_\epsilon=\epsilon  \mathbf e_1+\cdots +\epsilon \mathbf e_{m}+\lambda_1 \mathfrak a_1
+\cdots +\lambda_r \mathfrak a_{r},\\
 r\leq p+q-m, \epsilon>0, \lambda_i>0, all\ i.\end{array}
\end{equation}
be the varied vector as in (2.1). Also we denote a sub testing  map  
for the basis $$\mathfrak a_1, \cdots, \mathfrak a_{p+q-m}$$ by
$$D_{\bm\lambda}$$
with an arbitrary linear transformation. Now we choose the identity linear transformation for the basis
$$\mathbf e_1, \cdots, \mathbf e_{m}.$$ Then combine it with $D_{\bm\lambda}$ 
to obtain the testing  map 
$$D_{\bm\lambda_\epsilon}$$
for the linear space $V_k\times \mathbb R^{p+q-m}$ with
the center $$\mathbf u=u_1  \mathfrak a_1+\cdots+u_r\mathfrak a_{p+q-m}.$$
 Let $f(\mathbf v)d\mu_{p+q-m}$ be a test form on $ \mathbb R^{p+q-m}$. 
Then $$\psi =\theta_k(\mathbf v_k)\wedge f(\mathbf v)d\mu_{p+q-m}$$ is a test form 
on $V_k\times \mathbb R^{p+q-m}$. 
Then since $\mathcal L^k$  is of  Radon-Nikodym, 
$$\displaystyle{\lim_{|\bm\lambda_\epsilon| \Rsh 0}}\int_{V_k\times \mathbb R^{p+q-m}}
\mathcal L^k(D_{\bm\lambda_\epsilon}(\mathbf v_k, \mathbf v)) \psi(\mathbf v_k, \mathbf v)  d\mu_k\wedge d\mu_{p+q-m}
$$
exists.  Notice that  
\begin{align*}
& \displaystyle{\lim_{|\bm\lambda| \Rsh 0}}\int_{\mathbf v\in \mathbb R^{p+q-m}} \mathcal L^x_k(D_{\bm\lambda}(\mathbf v))
 f(\mathbf v)d\mu_{p+q-m}\\&
=\displaystyle{\lim_{|\bm\lambda_\epsilon |\Rsh 0}}\int_{V_k\times \mathbb R^{p+q-m}}
\mathcal L^k(D_{\bm\lambda_\epsilon}(\mathbf v_k, \mathbf v)) \psi(\mathbf v_k, \mathbf v)  d\mu_k\wedge d\mu_{p+q-m}.
\end{align*}
 where  $\mathcal L^x_k(\mathbf v)$ is the Radon-Nikodym derivative of 
$\mathcal I_k(\mathbf v)$ with respect to the Lebesgue measure, and $|\bm\lambda_\epsilon |\Rsh 0 $ is a particular zigzag path, 
\begin{equation}
 \displaystyle{\lim_{|\bm\lambda| \Rsh 0}}\  \displaystyle{\lim_{\epsilon\to  0}}.
\end{equation}

Thus  the $L^1$ function 
$\mathcal L^x_k(\mathbf v)$ is of  Radon-Nikodym. 
So is  the sum
$$\sum_{k} \mathcal L^x_k(\mathbf v).$$
We complete the proof.

\end{proof}

\bigskip

\begin{proposition} (intersection of the supports)
Let $T_1, T_2\in \mathcal C(\mathcal X)$. Then 
\begin{equation}
supp([T_1\wedge T_2])\subset supp(T_1)\cap supp(T_2).
\end{equation}

\end{proposition}

\bigskip

\begin{proof} Suppose
$$a\notin supp(T_1)\cap supp(T_2).$$
Then $a$ must be outside of either $supp(T_1)$ or $ supp(T_2)$. Let's assume first it is 
not in  $supp(T_2)$.  Since the support of a currents is closed, we 
choose a small  neighborhood $U_a$ of $a$ in  $\mathcal X$, but disjoint from $ supp(T_2)$.
Let $\phi$ be a $C^\infty$-form of $\mathcal X$ with a compact support in $U_a$. Then by Definition 3.1, part (2),  
when $\epsilon$ is small enough $R_{\epsilon} (T_2)$ is zero in $U_a$. Hence
\begin{equation}
\int_{ [T_1\wedge T_2]} \phi=0,
\end{equation}
for a test form $\phi$ supported in $U_a$.   Hence $a\notin  supp([T_1\wedge T_2])$.
 If $a\not\in supp(T_1)$, $U_a$ can be chosen disjoint with $supp(T_1)$. Then
 since $\phi\in \mathscr D(U_a)$ is a $C^\infty$-form of $\mathcal X$ with a compact support in $U_a$ disjoint with $supp(T_1)$, the restriction of $\phi$ to $T_1$ is zero. Hence 
$$
 \int_{ [T_1\wedge T_2]}\phi=0.
$$
Then $a\notin  supp([T_1\wedge T_2])$.
Thus $$a\notin supp(T_1)\cap supp(T_2)$$ will always imply 
$$a\notin  supp([T_1\wedge T_2]).$$

This completes the proof.

\end{proof}

\bigskip

\begin{ex}
Let $\mathcal X=\mathbb R^m$ be equipped with De Rham data consisting of single
open set with the convolution function $f$. Assume it has coordinates $x_1, \cdots, x_m$. 
   Let $$T_1=\delta_O dx_{1}\wedge \cdots\wedge dx_p, \quad 0<p<m$$ with the  $\delta$-function $\delta_O$
 at the origin $O$ of $\mathbb R^m$. Let $T_2$ be the $p$ dimensional plane $\{x_{p+1}=\cdots=x_{m}=0\}$. 
Now we consider the integral
\begin{equation}
\int_{T_1} R_\epsilon T_2.
\end{equation}
By the formula (3.5), it is equal to 
$$
 \int_{x\in T_1}\int_{y\in T_2=\mathbb R^{p}} {1\over \epsilon^m} f({x-y\over\epsilon})dx_{p+1}\wedge\cdots\wedge dx_{m}\wedge
dy_{1}\wedge\cdots\wedge dy_p.$$
By the continuity of the functional of the currents, we can interchange the order of $T_1, T_2$. Thus we first 
evaluate  $T_1$  at 
the differential form
$$ {1\over \epsilon^m} f({x-y\over\epsilon})dx_{p+1}\wedge\cdots\wedge dx_m$$
to obtain that 
\begin{equation}\begin{array}{c}
\int_{T_1} R_\epsilon T_2\\
\|\\
(-1)^{m(m-p) } \int_{y\in \mathbb R^{p}} {1\over \epsilon^{m}} f({-y_{1}\over\epsilon}, \cdots,{-y_p\over\epsilon}, 0, \cdots, 0 )
dy_{1}\wedge\cdots\wedge dy_p.
\end{array}\end{equation}
Since \begin{equation}\begin{array}{c}
\int_{y\in \mathbb R^{p}} {1\over \epsilon^{p}}f({-y_{1}\over\epsilon}, \cdots,{-y_p\over\epsilon}, 0, \cdots, 0 )
dy_{1}\wedge\cdots\wedge dy_p\\
=(-1)^p \int_{y\in \mathbb R^{p}} f(y_1, \cdots,y_p, 0, \cdots, 0 )
dy_{1}\wedge\cdots\wedge dy_p \end{array}\end{equation}
is a non-zero constant, $\int_{T_1} R_\epsilon T_2$ diverges to infinity as $\epsilon\to 0$.
This is consistent with Example 2.16 where $T_1$ does not satisfying the Lebesgue condition.
So a necessary condition for our convergence is the Lebesgue condition. 
\end{ex}

\bigskip

\begin{appendices}
\section{Kernel}

In [3]   De Rham created the notion of ``regularizing operator" which includes
De Rham's regulator  $R_\epsilon$.  Let $\mathcal X, \mathcal Y$ be two manifolds.  Let $L\in \mathscr D'(\mathcal X\times \mathcal Y)$.
There is a homomorphism
\begin{equation}\begin{array}{ccc}
\mathscr D(\mathcal X)\oplus \mathscr D(\mathcal Y) &\rightarrow & \mathbb R\\
(\phi_x, \phi_y) &\rightarrow &\int_{L}\phi_x\wedge \phi_y.\end{array}
\end{equation}
It leads to another homomorphism
\begin{equation}\begin{array}{ccc}
\Lambda; \mathscr D(\mathcal X)&\rightarrow &  \mathscr D'(\mathcal Y) \end{array}
\end{equation}
Then $L$ is called the kernel of $\Lambda$. Conversely given a homomorphism $\Lambda$, there is a kernel current $L$ on $\mathcal X\times \mathcal Y$. 
Notice $$\begin{array}{ccc}
\mathscr D(\mathcal X),  & &  \mathscr E(\mathcal Y)\\
\cap  & & \cap \\
\mathscr E'(\mathcal X), & &  \mathscr D'(\mathcal Y)
\end{array}$$
where $\mathscr E(\bullet)$ is the set of $C^\infty$ forms, and $'$ is the topological dual. 
\bigskip

\begin{definition}
(1) If $\Lambda$ can be extended to
  \begin{equation}\begin{array}{ccc}
\Lambda: \mathscr E'(\mathcal X)&\rightarrow &  \mathscr D'(\mathcal Y) \end{array}
\end{equation}
\hspace {1.5 CC} we say $\Lambda$ is regular.\bigskip

(2) If furthermore,  the regular $\Lambda$  has the image  inside of $\mathscr E(\mathcal Y)$, i.e.
\begin{equation}\begin{array}{ccc}
\Lambda: \mathscr E'(\mathcal X)&\rightarrow &  \mathscr E(\mathcal Y) \end{array}
\end{equation}
\hspace {1.5 CC}  we say $\Lambda$ is regularizing.

\end{definition}

\bigskip

\begin{theorem} (De Rham) 

$\Lambda$ is regularizing if and only if the kernel $L$ is a $C^\infty$ form. In particular
$R_\epsilon$ is regularizing. 
\end{theorem}

\end{appendices}

 \textsc{Department of Mathematics, Rhode Island college, Providence, 
   RI 02908}\par
  \text{E-mail address}:  \texttt{binwang64319@gmail.com}, \text{Fax}:   \texttt{1-401-456-4695}

\end{document}